\documentstyle[11pt]{article}
\def\df{\buildrel \rm def \over =}

\newtheorem{Theorem}{Theorem}
\newtheorem{Lemma}[Theorem]{Lemma}
\newtheorem{Proposition}[Theorem]{Proposition}
\newtheorem{Corollary}[Theorem]{Corollary}

\newfont{\Bbb}{msbm10 scaled\magstep1}	

        \newfont{\Bfit}{cmbxti10 scaled\magstep1}		
        
        \newfont{\Bfmit}{cmmib10 scaled\magstep1}		
        \newcommand{\bfmit  }[1]{\mbox{\Bfmit {#1}}}
        \newfont{\Sans}{cmss10 scaled\magstep1}			
        
        \newfont{\Frak}{eufm10 scaled\magstep1}			
        \newcommand{\frak}[1]{\mbox{\Frak {#1}}}

\newfont{\Bbbbb}{msbm7} 			
 
\newfont{\Bfmi}{cmmib7}

\def \Rk {{\Bbb R}^k}
\def \Rm {{\Bbb R} ^m}

\def \R{{\Bbb R}}
\def \Z{{\Bbb Z}}
\def \Zk{{\Bbb Z} ^k}
\def \Q{{\Bbb Q}}
\def \la{{\lambda}}
\def \Gm{{GL(m,{\Bbb R})}}

\def\df{\buildrel \rm def \over =}
\def\PROOF{{\em Proof}\,: }
\def\proof{{\em Proof}\,: }
\def\QED{~\hfill~ $\diamond$ \vspace{7mm}}

\def\Ad{\,\mbox{Ad}\,}

\def\ci{$C^{\infty}$}

\title{Differential Rigidity of Anosov Actions of Higher Rank Abelian Groups
and Algebraic Lattice Actions}
\author{A. Katok 
\thanks{Partially supported by  NSF grant DMS 9404061}
 and R. J. Spatzier 
\thanks{Partially supported by  NSF grant DMS 9626173}}

\date{}

\input amssymb.sty

\begin{document}
\maketitle

\vspace{.5em}
\begin{center}
{\em To Dmitry Viktorovich Anosov  on his sixtieth birthday}
\end{center}
\vspace{1em}

\begin{abstract}
We show that most homogeneous Anosov actions of higher rank Abelian groups
are locally $C^{\infty}$-rigid (up to an automorphism). This result is the main part
in the proof of local $C^{\infty}$-rigidity for two very different types of   algebraic actions of
irreducible   lattices in higher rank semisimple Lie groups: (i) the Anosov actions 
by automorphisms of tori and nil--manifolds, and (ii)   the actions of
cocompact lattices   on Furstenberg boundaries, in particular, projective spaces. The 
main new technical ingredient in the proofs is the use of a proper
``non--stationary'' generalization of the classical theory of normal forms for local 
contractions.
\end{abstract}

\section{Introduction}
The theory of Anosov systems is one of the crown jewels of dynamics.
The notion represents the most perfect kind of global hyperbolic
behavior which renders itself to qualitative analysis relatively easily
while the global classification remains to a large extent
mysterious. The concept was introduced by D.V.Anosov in \cite {A1}
and the  fundamentals of the theory were developed in his classic
monograph \cite {A2} which during almost thirty years that passed
after its publication has been serving
 as a source of  ideas and
inspiration for new work. Anosov himself called the class of systems
``U--systems'' after the first letter of the Russian word ``uslovije'',
which means simply ``condition'', and the name ``U--systems'' was being
used in the Russian language publications  for a number of years.
However, the term ``Anosov systems'', coined by Smale \cite {Smale},
who seems to have immediately recognized both the importance of the
notion and the credit due to its author, immediately became current in
publications outside of the Soviet Union and eventually was adopted
universally.

We are not going to discuss here either the history of the
theory of 
Anosov systems or its general impact in dynamics and beyond its
borders. We hope to  come back to these topics in another paper. We are
concerned here with only one line of development which has been
spawned by the seminal work of Anosov. The original notion of Anosov
system came in two varieties: Anosov diffeomorphisms (or ``cascades''
in Anosov's own terminology), i.e. the actions of the group $\Z$ of
integers, and Anosov flows, i.e. actions of the group $\R$. The former
are {\em structurally stable} as $\Z$ actions: every $C^1$
diffeomorphism  close to an Anosov diffeomorphism is {\em
topologically conjugate} to it. The latter allow countably many moduli
of topological conjugacy as $\R$ actions, e.g. the lengths of closed
orbits; they are nevertheless structurally stable in the orbit
equivalence sense: the orbit foliation of an Anosov flow is topologically
conjugate to the orbit foliation of any $C^1$ flow sufficiently close 
to it in the $C^1$ topology. The development  which concerns us here
deals with the notion of an {\em Anosov action} of a  group
(either discrete or continuous), more general than $\Z$ or $\R$. This
concept was originally introduced by Pugh and Shub \cite {PS} in the
early seventies (their notion coincides with ours  for abelian groups
and is different in general) 
and for many years fundamental differences between the ``classical''
($\Z$ and
$\R$) and other cases went mostly unnoticed, except for the
one--dimensinal situation \cite {Sacksteder}. Only with the progress
of the ``geometric rigidity'' approach to the actions of higher--rank
semisimple Lie groups and lattices in such groups initiated by Zimmer
in the 
mid-eighties
 did these differences come to the attention of
researchers beginning with \cite {KL1}. In particular, it turned out
that both local {\em differentiable} stability  (or rigidity) and
cocycle rigidity are present and in fact typical already for Anosov
actions of higher--rank abelian groups, i.e. $\Zk$ and
$\Rk$ for $k\ge 2$.
  
In this paper we bring the study of the local differentiable rigidity
of algebraic Anosov actions of $\Zk$ and $\Rk$ on compact
manifolds  as well as orbit foliations of such actions,  started in
\cite { KL1, KS1, KS2, KS4}, to a near conclusion.  Our results for the
abelian group actions formulated in Section 2.1 allow us  to obtain
comprehensive local
$C^{\infty}$ rigidity for two very different types of   algebraic
actions of irreducible   lattices in higher rank semisimple Lie
groups: (i) the Anosov  actions by automorphisms of tori and
nilmanifolds (Theorem 15, Section 3),  and (ii) the actions   on 
Furstenberg boundaries, in particular projective spaces (Theorem 17,
Section 4). While the latter area  was virtually unapproachable  save
for a  special result by Kanai
\cite {Kanai}, the former was extensively studied  beginning from  
\cite{Hurder1, KL1} and continuing in \cite{KLZ, Q1, Q4, Qian-Yue, KL2, Q2, Q3}.
Our result (Theorem 15 in Section 3) substantially extends all these
earlier works and brings the question to a  final solution. It also
gives further credence to the global  conjecture that all Anosov
actions of such lattices are smoothly equivalent to  one of the above.
Such a classification 
 was established for the much more special class of Cartan actions in
\cite{GS1}.

The objects considered in this paper, i.e. group actions
and foliations, are assumed  to be of class $C^ {\infty}$. 
Accordingly the basic notions, including the structural stability
are adapted to this case while they  sometimes appear more naturally
for the lower regularity.
In fact, for each of our results a certain
finite degree of regularity is sufficient  resulting in a loss of
regularity in the conjugating diffeomorphisms. We leave the detailed
study of optimal regularity conditions to a later paper.  In proper
places we make  specific comments about the degree of regularity
sufficient to guarantee existence  of $C^1$ conjugating
diffeomorphisms. 

Let $M$ be a compact manifold. We call two foliations ${\cal F}$ and ${\cal G}$ on
$M$ {\em orbit equivalent} if there is a homeomorphism $\phi : M \rightarrow M$
such that $\phi$ takes the leaves of ${\cal F}$ to those of ${\cal G}$. We  call
$\phi$ an {\em orbit equivalence}. If $\phi$ is a
$C^{\infty}$-diffeomorphism, we call ${\cal F}$ and ${\cal G}$
$C^{\infty}$-{\em orbit equivalent}. 

Endow the space  of foliations on $M$ with  the $C^1$-topology, i.e. two
foliations are close if their tangent distributions are $C^1$-close. We
call a
$C^{\infty}$-foliation ${\cal F}$ {\em structurally stable} if there is a neighborhood
$U$ of ${\cal F}$  such that any foliation in $U$ is orbit equivalent to ${\cal
F}$. We call a $C^{\infty}$-foliation ${\cal F}$ {\em $C^{\infty}$-locally rigid}
if there is a neighborhood $U$ of ${\cal F}$  such that any
$C^{\infty}$ foliation in $U$ is
$C^{\infty}$-orbit equivalent to ${\cal F}$.

Call a \ci -action  of a finitely generated discrete group $\Gamma$ 
{\em structurally stable} if any $C^{\infty}$ perturbation of
the        action  which is $C^1$-close on a finite generating set is 
conjugate to it via a homeomorphism. Call such a  $\Gamma$-action  
{\em locally $C^{\infty}$-rigid} if any  perturbation of the 
action  which is $C^1$-close on a finite generating set is 
conjugate to the original action  via a $C^{\infty}$ diffeomorphism.
We say that two actions of a group $G$ {\em agree  up to an automorphism}
if the second action can be obtained from the first by composition with an
automorphism of the underlying group. Call a \ci -action  of a Lie group $G$  
{\em locally $C^{\infty}$-rigid} if any  perturbation of the 
action  which is $C^1$-close on a compact generating set is 
$C^{\infty}$-conjugate up to an automorphism.  

The notions of $C^r,\,\,\,r\ge 1$ local rigidity for foliations
and actions are defined accordingly.   

For $\Gamma = \Z$ or $\R$, i.e. for the classical dynamical systems, 
diffeomorphisms
($\Z$-actions) and flows ($\R$-actions), local $C^{\infty}$-rigidity never 
takes
place. Moreover, the orbit foliations of flows are not locally 
$C^{\infty}$-rigid
either. In those cases, it does not help to allow the conjugacy in the 
definitions
to be only a $C^1$ diffeomorphism.  Structural stability on the other hand is a 
rather wide-spread
(although not completely understood) phenomenon. See \cite{KH} for a general
background and \cite{Robbin, Mane} for the definitive results on
structural stability in   $C^1$ (but
not $C^r,\,\, r\ge 2$)  topology. Since local
$C^{\infty}$-rigidity implies
structural stability  one can immediately see  from the necessary
conditions  for structural stability  that there are always moduli
of $C^1$ conjugacy  in the structurally stable case thus showing impossibility
of differentiable rigidity. 

Certain structural stability results, namely those  by Hirsch, Pugh and Shub
\cite{HPS}, are important for
our purposes. They  established  structural stability of central
foliations of certain partially hyperbolic dynamical systems which implies  
in particular that the orbit foliations of hyperbolic (Anosov) actions of $\Rk$
are structurally stable. At the beginning of Section 2.2 we explain how this fact
is used as the starting point in our proof of local differentiable rigidity.

Contrary to the classical case, in the higher--rank situation of a $\Z ^k$ or $\R
^k$-action for $k\ge2$,
$C^{\infty}$ local rigidity is possible and in fact seems to be as closely
related to the hyperbolic behavior as  structural stability is in the classical
case. Existence of this phenomenon was first discovered in
\cite [Theorem 4.2]{KL1} for certain (in a sense ``maximal'') actions of $\Zk$ by
toral automorphisms and was extended in \cite {KS4} to a  broader class of
Anosov actions  by both $\Rk$ and $\Zk$ (multiplicity--free standard actions).
The results on cocycle rigidity which first appeared in \cite {KS4} are
fairly definitive and they appeared in print as   \cite{KS1,
KS2}.  On the other hand, the local rigidity result for actions  in \cite {KS4}
was still too restrictive; in particularly, it covered the key situation of the Weyl
chamber flows only for split semi--simple groups. In the present paper we obtain
 much more general rigidity results which are close to  being definitive. We
prove 
$C^{\infty}$ local rigidity of most  known irreducible Anosov actions of  $\Zk$
and $\Rk$  as well as the orbit foliations of the latter. (Theorem 1, 
Corollaries 2 through 5). All such actions are essentially algebraic (see Section 2.1 
below)
and the remaining technical assumption of semisimplicity of the linear part of
the action is satisfied in  all interesting cases including the Weyl chamber flows,
which appear in applications to the rigidity of  actions of irreducible  cocompact
lattices in higher rank semisimple Lie groups on Furstenberg boundaries (Section 4),
and those actions by automorphisms of tori and nil--manifolds which are needed for
the proof of  rigidity of  actions of irreducible  
lattices in higher rank semisimple Lie groups by such automorphisms (Section 3).

We would like to thank E. Ghys, G. A. Margulis and G. Prasad for several 
helpful conversations.  A. Starkov 
pointed out to several inaccuracies in our original proof
of Corollary 4 and helped to correct those. We are also grateful to the Pennsylvania
State  University and the University of Michigan for their hospitality and  financial
support during several visits. Finally, the first author gladly acknowledges the support
from Erwin Schroedinger institute for mathematical physics where the final preparation of
this paper was completed.

The principal results of this paper are announced in \cite{KS5}.

\newpage

\section{Anosov actions of higher rank  Abelian groups}

\subsection{The main results}

We will consider ``essentially algebraic'' Anosov actions of either $\Rk$ or $\Zk$. 
We recall that an action of a group $G$ on  a compact manifold is {\em Anosov}
if some element $g \in G$ acts normally hyperbolically with respect to the 
orbit foliation (cf. \cite{KS1} for more details). To clarify
the notion of an  algebraic action, let us first define  affine algebraic
 actions of discrete groups.
Let $H$ be a connected Lie group with $\Lambda \subset H$ a cocompact lattice. 
Define $\mbox{Aff } (H)$ as the set of diffeomorphisms of $H$ which map right 
invariant vectorfields on $H$ to right invariant vectorfields. Define 
$\mbox{Aff } (H/\Lambda)$ to be the diffeomorphisms of $H/ \Lambda$ which lift
 to elements of $\mbox{Aff } (H)$. Finally, define an action $\rho$ of a 
 discrete group 
 $G$ on  $H/\Lambda$ to be {\em affine algebraic} if $\rho (g)$ is given by some
 homomorphism
$G \rightarrow \mbox{Aff } (H/ \Lambda)$. Let ${\frak h}$ be the Lie algebra 
of $H$. Identifying ${\frak h}$ with the right invariant vectorfields on 
$H$, any affine algebraic action determines a homomorphism $\sigma :
 G \rightarrow \mbox{Aut } {\frak h}$. Call $\sigma$ the {\em linear part }
of this action. We will also allow  quotient actions of these on finite 
quotients of $H/ \Lambda$, e.g. on infranilmanifolds. For  any Anosov 
algebraic action of a discrete group $G$, 
 $H$ has to be nilpotent (cf. eg. \cite[Proposition 3.13]{GS1}).

We define  algebraic $\Rk$-actions as follows. Suppose 
$\Rk \subset H$ is a subgroup of a connected Lie group $H$. Let $\Rk$ act on 
a compact quotient $H/\Lambda$ by left translations where $\Lambda$ is a 
lattice in $H$. Suppose $C$ is a compact subgroup of 
$H$ which commutes with $\Rk$. Then the $\Rk$-action on $H/ \Lambda$ descends 
to an action on $ C \setminus H / \Lambda$. The general algebraic 
$\Rk$-action $\rho$ is a finite factor of such an action. Let ${\frak c}$ 
be the Lie algebra of $C$. The {\em linear part} of $\rho$  is the 
representation of $\Rk$ on ${\frak c} \setminus {\frak h}$ induced by 
the adjoint representation
of $\Rk$ on the Lie algebra ${\frak h}$ of $H$.

Let us note that the suspension of an algebraic $\Zk$-action is an algebraic 
$\Rk$-action (cf. \cite[2.2]{KS1}).

The next theorem is our principal technical result for algebraic 
Anosov $\Rk$-actions.
 
We  denote the strong stable foliation 
of $a \in \Rk$ by ${\cal W} ^- _a$,  the strong stable distribution by $E^- _a$, and the
0-Lyapunov space by $E^0 _a$. Note that $E^0 _a$ is always integrable for
algebraic actions. Denote the corresponding foliation by ${\cal W} ^0 _a$. 
Also note that any algebraic action leaves a Haar measure $\mu$ on the quotient 
invariant. Given a collection of subspaces of a vectorspace, call a nontrivial 
intersection {\em maximal} if it does not contain any other nontrivial intersection of
these subspaces.

\begin{Theorem}   	\label{thm-main}
Let $\rho$ be an algebraic Anosov action of $\Rk$, 
for $ k \geq 2$, such that the linear part of $\rho$ is semisimple.
Assume that for any maximal nontrivial intersection 
$\bigcap _{i=1 \ldots r} {\cal W} ^- _{b_i}$
of stable manifolds of elements $b_1, \ldots, b_r \in \Rk$ there is an element 
$a \in \Rk$ such that for a.e. $x \in M$, 
$\bigcap _{i=1 \ldots r} E ^- _{b_i}  (x) \subset E^0 _a (x)$ and such that
a.e. leaf of the intersection $\bigcap _{i=1 \ldots r} {\cal W} ^- _{b_i}$ is
contained in an ergodic component of the one-parameter subgroup $t \, a$ of
$\Rk$ (w.r.t. Haar measure).

Then the orbit foliation of
   $\rho$ is locally $C^{\infty}$-rigid. In fact, the orbit equivalence
can be chosen $C^1$-close to the identity. Moreover, 
the orbit equivalence is transversally unique, i.e. for any two different 
orbit equivalences close to the identity, 
the induced maps on  the set of leaves  agree.
\end{Theorem}

The technical condition in this theorem guarantees that certain nontrivial 
intersections of stable manifolds are contained in the the closure of the 
orbit  of a generic point $x$ under one-parameter groups which translate 
these intersections by isometries.
It is similar to the condition  of the main
technical theorem of \cite{KS3}, Theorem 5.1, on invariant measures for
algebraic Anosov $\Rk$-actions.

\begin{Corollary}   	\label{cor-main}
Let $\rho$ be an algebraic Anosov action of $\Rk$, 
for $ k \geq 2$, such that the linear part of $\rho$ is semisimple.
Assume that  every one-parameter subgroup of $\Rk$ acts ergodically 
with respect to the Haar measure $\mu$. 
 Then the orbit foliation of
   $\rho$ is locally $C^{\infty}$-rigid. In fact, the orbit equivalence
can be chosen $C^1$-close to the identity. Moreover, 
the orbit equivalence is transversally unique.
\end{Corollary}

The assumption that any one-parameter subgroup of $\Rk$ acts ergodically is 
equivalent to the $\Rk$-action being weakly mixing, as one easily sees. 
We refer to Brezin and Moore for a more extensive discussion of the 
 ergodicity of homogeneous flows \cite{Brezin-Moore}. This condition 
is automatically satisfied if  $H$ is semisimple without compact factors
and $\Lambda$ is an irreducible lattice, as  follows
from Mautner's theorem \cite{Zimmer}. If $H$ is semisimple, then the linear
part of the action is automatically semisimple. The action of a split Cartan 
subgroup $A$ of $H$ on $K \setminus H /\Lambda$ where 
$\Lambda \subset H$ is a cocompact lattice in $H$, and $K$ is the compact part 
of the centralizer of $A$ in $H$ is always Anosov. This is the so-called {\em
Weyl chamber flow} \cite{KS1}. Hence the assumptions of
Corollary~\ref{cor-main} are
satisfied for the Weyl chamber flows.

More generally, we will establish \ci-local rigidity of  the twisted Weyl
chamber flows introduced in \cite{KS1}. Let us briefly recall their
construction. Let $A$ be a split Cartan subgroup of a real semisimple Lie group
$G$ without compact factors, and let  $K$ be the compact part of the
centralizer of $A$ in $G$. Let
$\Gamma$ be a cocompact irreducible lattice in $G$. Suppose that $\Gamma$ acts
on a nilmanifold $M =N/ \Lambda$ such that its linearization extends to a
homomorphism of $G$. Let $H$ be the semidirect product
$H= G\ltimes  N$. Then  $\Gamma\ltimes  \Lambda$ is a lattice in $H$. Then
$A$ acts on $K \setminus H / \Gamma\ltimes  \Lambda$. This action is Anosov
provided that the action of $\Gamma$ on $M$ contains an Anosov automorphism. In
this case, we call the $A$ action a {\em twisted Weyl chamber flow}.

Our conclusions for the Weyl chamber flows and twisted Weyl chamber flows may be
summarized as follows.

\begin{Corollary} 
The orbit foliation of the Weyl chamber flow on $K \setminus
H/\Lambda$, where $H$ is a semisimple Lie group of rank greater than
one without compact factors and 
$\Lambda \subset H$ is an irreducible lattice, is locally \ci -rigid.

The orbit foliation of a twisted Weyl chamber flows is locally  \ci
-rigid if the acting Cartan subgroup has real rank at least 2.

\end{Corollary}

In the  case  of Weyl chamber flows also any $C^6$ perturbation  of the orbit foliation
in $C^1$ topology is $C^1$ congugate to the  orbit foliation of the Weyl chamber flow.
For those groups where there are no positively proportional roots restricted to a split
Cartan (e.g. $SL(n,\R), n\ge 3$ or $SO(p,q),\,\, p,q\ge 2$) $C^6$
in this statement may be replaced by $C^4$. See Sections 2.2.1 and
2.2.3, Steps 4 and 5 for explanations.

To get rigidity results for algebraic $\Zk$-Anosov actions with
linear semisimple part,  we can   pass to the
suspension. Since the suspension of a $\Zk$ action is never weakly
mixing we cannot apply Corollary 2 and have to appeal to
Theorem~\ref{thm-main} in its full generality. As we will see, the
hypothesis  of this theorem is guaranteed if all nontrivial  elements
of the
$\Zk$-action are weak mixing. 
 By a theorem of W. Parry, an affine
 automorphism of a nilmanifold $H/\Lambda$ is weakly  
 mixing precisely when the 
quotient of the linear part on the abelianization 
 $\frac{\frak h}{[\frak h,\frak h]}$ does not have roots of unity as eigenvalues
\cite{Parry}.
As an orbit equivalence close to the identity between suspensions of 
$\Zk$-actions
  automatically  produces a conjugacy between the actions themselves,
we get the following corollary.

 \begin{Corollary}          \label{cor-nilautos}
 Let $\rho$ be an algebraic Anosov action of $\Zk$, $k \geq 2$,
  on an infranilmanifold with semisimple linear part. Suppose
    that no nontrivial element of the group   has  roots of unity
  as eigenvalues in the induced representation on the abelianization. 
   Then $\rho$ is \ci-locally rigid. Moreover, the conjugacy 
  can be chosen $C^1$-close to the 
identity, and is unique amongst  conjugacies close to the identity.
\end{Corollary}

All known Anosov $\Rk$-actions satisfying the assumptions of Theorem 1
belong to and almost exhaust the  list of {\em standard} $\Rk$-actions,
introduced  in
\cite{KS1}.  They  essentially consist of actions by
infranilmanifold (in particular, toral) automorphisms, Weyl chamber
flows, twisted  Weyl chamber flows and some further extensions (cf.
\cite{KS1} for more details ).

 For the standard Anosov 
actions, we showed that every smooth cocycle is smoothly cohomologous 
to a constant cocycle \cite[Theorem 2.9]{KS1}. As a consequence, all 
smooth time changes are smoothly conjugate to the original action
 (possibly composed with an automorphism of $\Rk$). Combining this 
 with the above results about the \ci rigidity of the orbit foliations, we obtain the
following corollary.

\begin{Corollary}
The standard algebraic Anosov actions of $\Rk$ for $k \geq 2$ with 
semisimple linear part are  locally $C^{\infty}$-rigid. 
Moreover, the $C^{\infty}$-conjugacy $\phi$ between the action composed with 
an automorphism $\rho$ and a perturbation can be chosen $C^1$-close to the 
identity.
The automorphism $\rho$ is unique and also close to the identity. Finally,
 $\phi$  is  unique amongst  conjugacies close to the identity modulo
  translations in the acting group. 
\end{Corollary}

In the case of Weyl chamber flows, there are in fact very few conjugacies, 
close to the identity or not. We pursue this in Section 2.5.

For  our application  to the local rigidity of projective lattice
actions, we will  need a rigidity result similar to Corollary 3
 for Anosov actions of certain reductive groups.
 Let $G$ be a connected semisimple Lie group with  finite center and 
without
compact factors. Let $\Gamma \subset G$  be an irreducible cocompact lattice. 
Let $P$ be 
a parabolic subgroup of $G$, and let $H$ be its Levi subgroup. Thus $P= H \,
U^+$ where $U^+$ is the unipotent radical of $P$. 
(We refer to \cite{Warner} for details on parabolic subgroups and boundaries
of $G$.) Then $H$ acts on $G/\Gamma$ by left translations. These actions are
Anosov (cf. eg. \cite{Leuzinger}).

\begin{Theorem} If the real rank of $G$ is at least 2, then the orbit foliation
${\cal O}$ of $H$ is $C^{\infty}$-locally rigid. Moreover, the orbit equivalence
can be chosen $C^1$-close to the identity.
\end{Theorem}

 We will actually use the following corollary of the proof of the  theorem. 
 Let $P=L\, C\, U^+$ be the Langlands  decomposition of $P$ (with respect to some Iwasawa
 decomposition $G=K\, A\, N$ of $G$), and let $M_P$ be the
 centralizer of $C$ in a  $K$. Then the orbit foliation
of $H$ on $\Gamma \setminus G $ descends to a foliation ${\cal R}$ on 
$\Gamma \setminus G /M_P$. 

\begin{Corollary} If the real rank of $G$ is at least 2, then the 
foliation ${\cal R}$ on $\Gamma \setminus G /M_P$ is $C^{\infty}$-locally rigid.
Moreover, the orbit equivalence
can be chosen $C^1$-close to the identity.
\end{Corollary}

\subsection{Proof of Theorem 1}

 The general outline of the proof is  as follows. The weak stable
foliations of various elements of our action are canonically perturbed to foliations with smooth
leaves saturated by the perturbed foliation ${\cal R}$. Due to the
Hirsch-Pugh-Shub structural stability theorem for normally hyperbolic
systems, there is  a H\"{o}lder orbit equivalence between
${\cal O}$ and ${\cal R}$. This orbit equivalence is smooth along
${\cal O}$, and carries the weak stable foliations of various elements
into their canonical perturbations. By transversality, the same is
true for the intersections of the weak stable foliations. The tangent
bundles of suitable intersections in fact provide a transversal
splitting of the tangent bundle.  The arguments outlined so far are
fairly standard. In particular, they form a basis of our earlier
results
\cite{KS4} including the proof of local differential rigidity in the multiplicity free case  as well
as   partial results in more general cases.

The key part in the proof is the smoothness of the orbit equivalence along these intersections.
The essential new  technique used here  is a suitable non-stationary generalization of the classical
theory of normal forms for local differentiable contractions \cite{K2}.  This technique is
summarized in the following section. Once  we have smoothness of the orbit equivalence along  the
foliations of the splitting, standard elliptic theory shows smoothness.
Theorem 6 dealing with a more general reductive case is proved very
similarly so there is no need to repeat all the details. We outline
the necessary modifications after the proof of Theorem 1.

\subsubsection{Preliminaries on normal forms}

We now give a summary of one possible non-stationary generalization of
the  classical normal form theory for local contractions whose origins
go back to Poincar\`{e} and which, in its modern form for the
$C^{\infty}$ case, was developed by Sternberg and Chen
\cite{Sternberg,Chen}. It is quite possible that one can find results
very similar to the ones below in the vast literature on normal
forms.  The first author wrote the precise versions needed for our
application in
\cite{K2}.

Consider a continuous extension ${\cal F}$ of a homeomorphism $f$ on a compact connected metric
space $X$ to  a neighborhood of the zero section of a  vector bundle $V$ which is smooth along the
fibers and preserves the zero section. Let $F = D {\cal F} _0$ where
the derivative is taken at the zero section in the fiber direction.
Consider the induced operator $F ^*$ on the Banach space of continuous
sections of  the bundle $V$ endowed with the uniform norm given by 
$F^* \, v (x) =  F (v (f^{-1} (x)))$. 
 The {\em characteristic set} of $F ^*$ is the set of logarithms of the absolute values of the
spectrum of $F ^*$. It consists of the union of finitely many intervals. 
For notational convenience we will consider a slightly more general situation. Namely, 
consider a finite set of disjoint intervals 
 $\Delta _i = [ \lambda _i, \mu _i ]$  on the negative half--line  which {\em contain} 
 the characteristic set.  We order the intervals in increasing order such that
$\lambda _{i+1} > \mu _i$.    Then the bundle $V$ splits into the direct sum of
$F$--invariant subbundles $V_1,\dots,V_l$  such that the spectrum of the 
restriction 
$ F\mid _{V_i} \subset \Delta_i,\,\,\,\,i=1,\dots,l$ \cite{K3}. Let $m_i$ be the
 dimension of the sub-bundle
$V_i$.
 Call the extension ${\cal F}$ a {\em contraction} if $F=D {\cal F}_0$ is a contraction with respect
to a continuous family of Riemannian metrics in the fibers. This is equivalent to the condition that 
the $\log$ of the
spectral radius of $F^*$ is negative, i.e. $\mu _l <0$.

 We will assume that  $F ^*$ has {\em narrow band spectrum}, i.e. that
\[ \mu _i + \mu _l < \lambda _i \] for all $i=1,\ldots,l$ \cite{K2}.   

Represent the space $\Rm$ as the direct sum of subspaces $V_1,\dots,V_l$ of dimension 
$m_1,\dots, m_l$ correspondingly and 
let $(t_1,\dots,t_l)$ be the
corresponding coordinate representation of a vector $t\in\Rm$.
 Let $P:
\Rm\to\Rm;\,\,(t_1,\dots,t_l)\mapsto (P_1(t_1,\dots,t_l),
\dots,P_l(t_1,\dots,t_l))$ be a polynomial
map preserving the origin. We will say that the map $P$ is of {\it  
sub--resonance type} if it has  non--zero homogeneous terms in
$P_i(t_1,\dots,t_l)$ with degree of homogenuity $s_j$ in the coordinates of
$t_j,\,\,i=1,\dots,l$ only if 
$$ \hspace{15em}   \la_i\le\sum s_j\mu_j \hspace{15em}   (*)$$

Let us notice that the notions of a homogeneous polynomial and degree of homogenuity in 
a Euclidean space are invariant under linear transformations. Thus the notion of a map of
sub--resonance type depends only on the decomposition $\Rm=\bigoplus_{i=1}^lV_i$,
but not on a choice of coordinates in each component of this decomposition. 
Furthermore, if
$V$ is a continuous vector bundle which splits into  the direct sum of continous
subbundles $V_1,\dots,V_m$, and if the vectors $\la$ and
$\mu$ are given one can  define in a natural way the notion of a continuous polynomial bundle
 map of sub--resonance type. 

 We will call any
 inequality  of   type $(*)$ a {\it  sub--resonance relation}. There are always
sub--resonance relations of the form
$\la_i\le\mu_j$ for $j=i,\dots,l$. They correspond to the  linear terms of the polynomial.
We will call such relations {\it trivial}.  The narrow band condition guarantees that for
any non--trivial sub--resonance relation $s_j=0$ for $j=1,\dots,i$.

Although the composition of maps of sub--resonance type may not be a map of sub--resonance 
type  it is not difficult to see that maps of sub--resonance type with invertible derivative 
at the origin generate a finite--dimensional group all of whose elements are polynomial maps 
of degree bounded by a constant which depends only on the vectors $\la$ and $\mu$.
We  will denote this group by
$SR_{\la,\mu}$ and will call its elements {\it sub--resonance generated} polynomial maps. 
See \cite{K2} for more explanations.

In certain cases the maps of sub--resonance type with invertible derivative at the origin 
already form a group. In particular, if the whole spectrum is sufficiently narrow i.e if
$\la_1>2\mu_l$, then there are no non--trivial sub--resonance relations and hence 
$SR_{\la,\mu}$ is a subgroup of $\Gm$ the group of linear automorphisms of
$\Rm$. The next simplest case appears in the perturbation of the point spectrum with
$2:1$ resonance. In this case
$l=2$ and  $\la_1>\max \{3\mu_2, \mu_1+\mu_2\}$. The only possible non--trivial
sub--resonance relation is $\la_1\le 2\mu_2$. In this case
 the group
$SR_{\la,\mu}$ consist of quadratic maps of the form
$$
 P(t_1,t_2)=(L_1 t_1 +Q(t_2,t_2),L_2t_2),
$$  where $L_1$ and $L_2$ are linear maps and $Q$ is a quadratic form. The above two cases
are the only ones which appear in the consideration of small perturbations of  Weyl
chamber flows since only the  double of a root (restricted to a split Cartan) can be a (restricted) root.

Call two extensions {\em conjugate}  if there exists a continuous family of local 
$C^{\infty}$ diffeomorphisms of the fibers $V(x)$, preserving the origin which transforms
one extension into the other. The following two theorems on normal forms and centralizers
are  proved in \cite{K2}.

\begin{Theorem} {\bf (sub--resonance normal form)} Suppose the extension ${\cal F}$  is a
contraction and the linear extension 
$D\cal F_0$ has a narrow band spectrum determined by the vectors $\la=(\la_1,\dots,\la_l)$
and
$\mu=(\mu_1,\dots,\mu_l)$. There exists an extension $\tilde{\cal F}$ 
conjugate to $\cal
F$  such that for every $x\in X$,
$$
\tilde{\cal F}\mid _ {V(x)} :\,\,\, \bigoplus_{i=1}^lV_i(x)\to \bigoplus_{i=1}^lV_i(f(x))
$$ is given by a polynomial map of sub--resonance type. 
\end{Theorem}

\begin{Theorem} {\bf (Centralizer for sub--resonance maps)} Suppose $g$ is a homeomorphism
of the  space
$X$ commuting with $f$ and $\tilde{\cal G}$ is an extension of $g$ by
$C^{\infty}$--diffeomorpisms of the fibers (not necessarily a contraction) commuting with
an extension $\tilde{\cal F}$  satisfying the assertion of Theorem 8. Then $\tilde{\cal
G}$ has a similar form:
$$
\tilde{\cal G}\mid _ {V(x)} :\,\,\,\bigoplus_{i=1}^lV_i(x)\to \bigoplus_{i=1}^lV_i(g(x))
$$ is a polynomial  sub--resonance generated map, i.e. an element of the group
$SR_{\la,\mu}$. 
\end{Theorem}

Combining these two theorems we see that a local action of an abelian group by extensions
which contains a contraction which has  narrow
 band spectrum can be simultaneously brought to a normal form. 

\begin{Corollary}        \label{cor-jointnormal} Let $\rho$ be a continuous action of
$\Rk$ on a compact connected metric space $X$. Let $V$ be a vector bundle over $X$.
Suppose that $\sigma$ is a local action of  $\Rk$ in  a neighborhood of the zero section
of $V$ such that $\sigma$ covers
$\rho$, $\sigma$ is \ci along the leaves and each 
$\sigma (a) \mid _{V_x}$ depends continuously on the base point $x$ in the
$C^{\infty}$-topology. Suppose further that some $a \in \Rk$, $\sigma (a)$ is a
contraction and that the induced linear operator on $C^0$-sections of $V$ has narrow band
spectrum.   Then there exist  $C^{\infty}$ changes of coordinates in the fibers $V_x$,
depending continuously on $x$ such that for all $b \in
\Rk$, $\sigma (b) \mid _{V(x)}$ is a polynomial map of sub--resonance type.
\end{Corollary}

All these results have counterparts for extensions of finite differentiability
(M.Guysinsky; in preparation). The precise degree of differentiability  both for the
normal form and for the centralizer  theorems  depends on the vectors $\lambda$ and $\mu$.
 Let us mention the following specific cases which appear  in the perturbations of the
Weyl chamber flows  (see comments after the formulation of Corollary 3)
 and hence in the applications to the projective lattice actions on the Furstenberg
boundaries (See Section 4). We will refer to  the differentiability of an extension
meaning a leafwise differentiability continuously  changing in the corresponding topology
and to differentiability of a normal form  meaning  the leafwise differentiability of the
map bringing an extension to the normal form.

In the narrow spectrum case
$\lambda_1>2\mu_l$ there is a $C^2$ linear normal form  for a 
$C^4$ extension and the   $C^2$ centralizer of a linear extension with this kind of
spectrum consist of linear maps.  In the case $l=2$ and 
 $\la_1>\max \{3\mu_2, \mu_1+\mu_2\}$ there is a $C^3$ sub--resonance  normal 
 form for a $C^6$ extension and any commuting $C^3$ extension has the same 
 sub--resonance form.

\subsubsection{Coarse Lyapunov decomposition for algebraic actions}

Consider an algebraic $\Rk$-action $\rho$ with  linear part $\sigma$.  We will assume that
$\sigma$ is semisimple. Define the {\em Lyapunov exponents} of $\rho$ as the $\log$'s of
the absolute values of the eigenvalues of $\sigma$. We get linear functionals $\chi : \Rk
\rightarrow \R$. These functionals coincide with the Lyapunov exponents for the Haar
measure
\cite{KS1,K2} or for any invariant measure for that matter. Note that
$\rho$ is Anosov precisely if the only 0 Lyapunov exponents come from the orbit.  There is
a splitting of the tangent bundle into $\Rk$-invariant subbundles 
$TM = \oplus _{\chi} E_{\chi}$ such that the Lyapunov exponent of $v \in E_{\chi}$ with
respect to
$\rho (a)$ is given by $\chi (a)$. We call $E_{\chi}$ a {\em Lyapunov space} or {\em
Lyapunov distribution} for the action. Then the strong stable  distribution $E^- _a$ of $a
\in \Rk$ is given by $E^- _a = \sum _{\chi (a) <0}  E_{\chi}$. While individual Lyapunov
distributions may not be integrable, the sum 
$E^{\chi} =\oplus E_{\lambda}$ where $\lambda$ ranges over all Lyapunov functionals which
 are positive multiples of a given Lyapunov functional $\chi$ is always integrable. In
fact, set
$H=\{a \in \Rk \mid \chi (a) \leq 0\}$, and call it a {\em Lyapunov half space}. Then $E^{
\chi}$ is the intersection of all stable distributions of elements in $H$, as is easily
seen. We will also  denote $E^{ \chi}$ by $E_H$.  This is an integrable  distribution with
integral foliation ${\cal W} ^- _H$ whose leaves are intersections of strong stable
manifolds. Note that
${\cal W} ^- _H$ is also integrable with the orbits. Denote the resulting foliation by 
${\cal W} _H$. Then we get a decomposition 
\[ TM = \bigoplus E_H  \oplus T{\cal O}\] where $H$ runs over all Lyapunov half spaces and
$T {\cal O}$ is the tangent bundle to the orbits of the action. We call this
decomposition  the {\em coarse Lyapunov decomposition} of $TM$. Note that each $E_H$ is
the sum of Lyapunov distributions corresponding to Lyapunov exponents proportional to each
other with positive coefficients of proportionality.  Thus for any Lyapunov half space $H$
one can find a uniquely defined Lyapunov characteristic exponent $\chi(H)$ (called the
{\it bottom exponent} for $H$) and positive numbers
$1=c_1(H)<c_2(H)<,\dots,<c_{m(H)}(H)$ such that 
$$ E_{H}(x)=\bigoplus_{i=1}^{m(H)}E_{c_i(H)\chi(H)}.
$$ The situation is particularly simple for algebraic $\Rk$-actions on  quotients of
semisimple groups, i.e. for Weyl chamber flows. As was noticed in the previous subsection
the only positive coefficient proportionality with a bottom exponent in that case is two.
This  fact is important for the case of low smoothness mentioned above.

Finally note that the coarse Lyapunov spaces are precisely the maximal nontrivial 
intersections of  stable distributions of elements in $\Rk$.

\subsubsection{Existence and regularity of the orbit equivalence}
\vspace{1em} 

{\bf Step 1: H\"{o}lder conjugacy.} Now consider an algebraic Anosov action $\rho$ of
$\Rk$ on a manifold $M$. Let ${\cal O}$ denote its orbit foliation. Suppose ${\cal R}$ is 
a  $C^1$-small perturbation of ${\cal O}$. Let $a$ be a normally hyperbolic
 element of $\Rk$. Then the 1-parameter group $t\, a $ defines a vectorfield 
$V$ on $M$. Define a perturbation $\tilde{V} (p)$ for any point $p \in M$ by projecting
$V(p)$ to $T {\cal R} _p$, the tangent space to ${\cal R}$ at $p$.  By Hirsch, Pugh and
Shub's structural stability theorem,  $\tilde{V}$ is  normally hyperbolic to a   foliation
$\tilde{{\cal R}}$, invariant under 
$\tilde{V}$ which is $C^1$-close to 
${\cal O}$ \cite[Theorem 7.1]{HPS}. Since $\tilde{V}$-invariant foliations $C^1$-close to 
${\cal O}$ are unique by normal hyperbolicity, ${\cal R}=\tilde{{\cal R}}$, and thus
$\tilde{V}$ is normally  hyperbolic to ${\cal R}$. Again applying structural stability, 
there exists a H\"{o}lder orbit equivalence $\phi$ between ${\cal R}$ and ${\cal O}$. We
may also assume that
$\phi $ is
$C^{\infty}$ along the leaves. In fact,  one identifies the image of a leaf under $\phi$
uniquely by shadowing (this also proves the transversal uniqueness of the orbit
equivalence). Then one  can map a leaf to its image by intersecting the image of a
neighborhood of the  zero section in the normal bundle to the leaf with the image leaf.
This yields an orbit equivalence smooth along the leaves. 

Let us also note here that to prove smoothness of $\phi$ it suffices to prove transversal
smoothness of $\phi$. More precisely, pick smooth  transversals $T_1$  and $T_2$ to ${\cal
O}$ and ${\cal R}$ respectively. Map $T_1$ to $T_2$ by sending 
$t \in T_1$ to the intersection of the leaf of ${\cal R}$ through $\phi (t)$  with $T_2$.
Call
$\phi$ {\em transversally smooth} if the resulting map from $T_1 \rightarrow T_2$ is
$C^{\infty}$. Now the above construction shows that transversal smoothness of $\phi$
implies that $\phi$ is $C^{\infty}$.
\vspace{1em} 

{\bf Step 2: Invariance of weak stable foliations.} Now consider the coarse Lyapunov
decomposition of $\rho$. We will first show that 
$\phi$ takes the foliations ${\cal W} _H$ into  $C^1$-close foliations $\tilde{{\cal W}}
_H$ saturated by ${\cal R}$.  Pick normally hyperbolic elements $a_1,\ldots , a_q$ such
that $E^-_H$ is the intersection of the strong  stable distributions $E^- _{a_s},\,
s=1,\ldots, q$. Let $V_s$ denote
 the vectorfields that generate the flows $\rho(t \, a_s)$ on $M$. As above, let 
$\tilde{V} _s$ be the projection of $V_s$ to the tangent distributions of ${\cal R}$ 
(Note that the
$\tilde{V} _s$ do not necessarily commute). Nevertheless, ${\cal R}$  is invariant under
$\tilde{V} _s$ for all $s$, and all $\tilde{V} _s$ are normally hyperbolic with respect to
${\cal R}$. As $\tilde{V} _s$ and $V_s$ are $C^1$-close, structural stability implies that
$\phi$ maps the weak stable foliations of $V_s$ to
 those of $\tilde{V} _s$ \cite{HPS}. As the tangent distributions of  the weak stable
foliations of  $\tilde{V} _s$ are 
$C^0$-close to those of $V_s$, they still intersect transversely in a foliation
$\tilde{{\cal W}} _H$ which by definition is saturated under ${\cal R}$. In particular, the
$\tilde{{\cal W}} _H$ are invariant under the flows of 
 $\tilde{V} _s$, and in particular under their time 1 maps $\tilde{a} _s$.  Denote their
tangent distributions by $\tilde{E} _H$.

Consider the induced operator $a_s ^*$ on the Banach space of sections of  the bundle $E_H$
given by
$a_s ^* \, v (x) = D\rho( a_s) (v \rho(a_s ^{-1} (x))$.  Then the characteristic set of
$a_s ^*$ just consists of the numbers $\exp \chi (a_s)$ where $\chi$ runs over all the
Lyapunov exponents with
$E_{\chi}\subset E_H$.
\vspace{1em} 

{\bf Step 3: Construction of a perturbed action and natural
extension.} Define a continuous action
$\tilde{\rho}$ of $\Rk$ on $M$ by conjugating $\rho$ by $\phi$. Next
we will construct an extension of this action on the  bundle $T$,
defined as follows. Since $\tilde{{\cal W}}_H (x)$ is $C^1$-close to
${\cal W} _H  (\phi ^{-1} (x))$, we can project ${\cal W} _H (\phi
^{-1} (x))$ to 
$\tilde{{\cal W}}_H (x)$ by a $C^{\infty}$ map, close to the
identity.  In particular, the image of the foliation ${\cal W}^- _H$
on 
${\cal W} _H (\phi ^{-1} (x))$ defines a smooth foliation ${\cal T}$
with leaves
$T_y, \, y \in \tilde{{\cal W}}_H (x)$, such that the leaf $T_x$ is
$C^0$-close to $\phi ({\cal W} ^- _H (\phi ^{-1} (x))$. Locally we can
identify the bundle of $T_x$ with the normal bundle of 
${\cal R}$ in $\tilde{{\cal W}}_H $.  Finally, define the extension of
$\tilde{\rho}$ by holonomy as follows. If $y \in T_x$ and $a \in \Rk$,
let ${\cal A} (y)$ be the unique  intersection point of the local leaf
of ${\cal R} (\tilde{\rho}(a) (y))$ with $T_{\tilde{\rho} (a) (x)}$.

Let ${\cal A} _s$ denote the extensions of  $\tilde{\rho}(a_s)$ on
$T$. These extensions commute with each other. Define the operators
${\cal A} _s^*$ on the Banach space of continuous sections of $T$,
endowed with the uniform norm, as usual. Since the extensions ${\cal
A} _s$ are  small perturbations of the
$\rho(a_s)$, the spectra of the ${\cal A} _s^*$  are close to those of
$a_s ^*$. Thus  the characteristic set of  ${\cal A} _s^*$  is 
contained in small intervals $\Delta _i = [ \lambda _i, \mu _i ]$
about the characteristic set of
$a_s ^*$ where we order the intervals in increasing order such that
$\lambda _{i+1} > \mu _i$.   The size of the intervals depends on the
size of the perturbation of $\rho$. Thus we may assume that the 
operator
${\cal A} _s^*$ has narrow band spectrum (cf. 2.2.1).
\vspace{1em} 

{\bf Step 4: Smoothness along coarse Lyapunov directions - the main
step.} We will now show that
$\phi: {\cal W} ^- _H \rightarrow \tilde{{\cal W}} _H $ is smooth. 
 Let $\psi _x : {\cal W} ^- _H (x) \rightarrow T_{\phi (x)} $ denote
the composition of $\phi$ and projection to $T_{\phi (x)}$ along the
leaves of ${\cal R}$.  

As  we will see, these $\psi _x$ intertwine the smooth transitive
actions of a Lie group $G$ on both
$ {\cal W} ^- _H (x)$ and $T_x$. Then conjugation by $\psi _x$ defines
a continuous and hence
$C^{\infty}$ homomorphism $\eta _x : G \rightarrow G$. Since 
$\psi _x (g \, x) = \eta _x (g) \, \psi _x (x)$, we see immediately
that $\psi _x$ is 
$C^{\infty}$.

Let us first define $G$ and its action on  $ {\cal W} ^- _H (x)$  for
a.e. 
$x \in M$.  Let $\chi$ be a  Lyapunov exponent such that $H$ is the
halfspace on which $\chi$ is
 nonpositive.  Since $ {\cal W} ^- _H (x)$ is a maximal nontrivial 
  intersection of stable
 manifolds, by the hypothesis of the theorem there is a one-parameter
subgroup $t\, b$ such that $\chi (b) =0$ and whose ergodic components
are saturated by  leaves of $ {\cal W} ^- _H $.
 Let $M$ be the quotient of a compact homogeneous $H/\Lambda$ by a
compact 
 group $K$ commuting with an affine action $\rho ^*$ of $\Rk$ such
that $\rho$ is the quotient action.  Then we can endow $H$ with a
right invariant metric such that $\rho ^* $ dilates this metric 
according to the linear part
$\sigma ^*$ of $\rho ^*$ and $K$ leaves the metric invariant. Let
$\bfmit g$  be the induced Riemannian metric on $M$.
 Then the maps $\rho(t b): {\cal W} ^- _H (x) \rightarrow {\cal W} ^-
_H (\rho(t b(x)))$ are isometries with respect to the  Riemannian
metrics induced by  $\bfmit g$ on the 
$ {\cal W} ^- _H (x)$. 

Since the ergodic components of $\rho(t b)$ are saturated by leaves of 
$ {\cal W} ^- _H $,    for a.e. point $x$, the closure of the orbit
$\rho(t b)(x)$ contains
 $ {\cal W} ^- _H (x)$. Fix such a point $x$. In particular,  given
any $ y \in  {\cal W} ^- _H (x)$, there is a sequence of $t_n$ such
that 
$\lim _{n\rightarrow
 \infty}\rho( t_n b ) (x) =y$. We may assume moreover, that the
sequence of
  isometries
$\rho( t_n b ) : {\cal W} ^- _H (x) \rightarrow {\cal W} ^- _H\rho(
t_n b )(x) )$  converges to an isometry $g: {\cal W} ^- _H (x)
\rightarrow {\cal W} ^- _H (x)$ which takes
 $x$ to $y$. In fact, it is easy to see that such a limit is a local
isometry
 and also a covering map.
 Since ${\cal W} ^- _H (x)$ is simply connected, it follows that the
limit is a
 global isometry (cf. \cite[Proposition 2.9]{GS1} for a  detailed
proof). 
 Let $G = G _x $ be  the closure of the  group generated by all such
limits $g$ in the isometry group of ${\cal W} ^- _H (x)$. Then $G$
acts transitively on $ {\cal W} ^- _H (x)$.  

 Since $\psi _x : {\cal W} ^- _H (x) \rightarrow T_x$ is a
 homeomorphism,  $G$ acts on $T_x$ by
conjugating by $\psi _x$. We will show that for all  $g \in G$,
 $g$ acts
 smoothly on $T_x$.  
Let us first show this for a limit $g$
of  isometries $g_n = \rho( t_n b ) : {\cal W} ^- _H (x) \rightarrow 
{\cal W} ^- _H (\rho( t_n b ) (x) )$. Then $h_n := \psi _{g_n \, x} 
\, g_n \, \psi _x ^{-1}$ is a compact family of maps in the 
$C^0$-topology, converging to a map $h: T_{\phi (x) }\rightarrow  T_{\phi (x)}$.
 Note that $h$ is
$g$ acting on $T_{\phi (x)}
$. On the  other hand, the $h_n$ are $C^{\infty}$-maps, 
as in fact they are holonomy maps for ${\cal R}$. 

Now comes the central point in the argument. Since natural
extension 
${\cal A}$ of the conjugated action $\tilde{\rho}$ contains  a
contraction ${\cal A}_s$ satisfying the narrow band condition, by
Corollary~\ref{cor-jointnormal} we can find a family of
$C^{\infty}$ coordinate changes in the fibers $T_x$, depending
continuously on $x$, such that the 
${\cal A}_s$, and therefore the $h_n$, which commutes with ${\cal
A}_s$, written in these coordinates act as elements of a certain finite
dimensional Lie group, namely the group $SR_{\lambda,\mu}$. As the
$h_n$ are compact in
$C^0$, they are also a compact family in $SR_{\lambda,\mu}$. Therefore,
the
$h_n$ converge in the
$C^{\infty}$-topology, and $h$ is smooth, and belongs to 
$SR_{\lambda,\mu}$. Moreover, the group generated by such $h$ as well
as its closure $G$ belongs to $SR_{\lambda,\mu}$. Thus every element in
$G$ acts smoothly on $T_x$. 

Now suppose that $y$ is a limit of points $x_n$ such that ${\cal W} ^-
_H (x_n)$ is contained the orbit closure of $x_n$ under $\rho(t b)$.
Then the special Lie 
 group of isometries $G_{x_n}$ constructed above is transitive on
${\cal W} ^- _H (x_n)$, and acts smoothly on $T_{x_n}$ by
diffeomorphisms in $SR_{\la,\mu}$.  Furthermore, $\psi _{x_n}$
conjugates these two actions. Let $G_y$ be the group generated by all
$C^0$-limits of $g_n \in G_{x_n}$. Then $G_y$ is a transitive group of
isometries of ${\cal W} ^- _H (y)$. Let $g = \lim g_n$ be such a
limit.  Then $\psi _{x_n} \circ g_n \circ \psi _{x_n} ^{-1}$ converges
to $\psi _y \circ g \circ \psi _y ^{-1}$ in the $C^0$-topology. Since
the 
$\psi _{x_n} \circ g_n \circ \psi _{x_n} ^{-1}$ again belong to
$G_{\lambda,\mu}$, this convergence is again \ci. Finally, let us note
that any
$y \in M$ is such a limit since the union of the ergodic components of
$\rho(t b)$ have full measure. 

Thus  we have shown that for all $y \in M$, there is a Lie group 
$G_y $ of isometries of ${\cal W} ^- _H (y)$ which acts continuously
on the  transversals $T_y$  by diffeomorphisms and such that $\psi _y$
conjugates these actions. By a theorem of Montgomery, the $G_y$ action
on the 
$T_y$ is  $C^{\infty}$ \cite[Section 5.1, Corollary]{Mont-Zipp}. As we 
mentioned  at the beginning of this step of the proof, this is
sufficient to conclude that  $\psi_y$ is a 
$C^{\infty}$ diffeomorphism since  $\psi_y$  intertwines two smooth
transitive actions of a Lie group.

In the case of $C^6$  (of $C^4$ in the absence of double roots)
perturbations of Weyl chamber flows we use  instead of Corollary 10 the
$C^3$ (corr.
$C^2$) normal form and the rigidity of the $C^3$ (corr. $C^2$) centralizer
established by Guysinsky (see the end of Section 2.2.1). The
conclusion is that the conjugation intertwines two transitive Lie groups
of diffeomorphisms of corresponding finite regularirty. The finite
regularity version of the Montgomery theorem implies then that
the conjugacy is at least $C^1$ along $\cal W^-_H$.

\vspace{1em}

{\bf Step 5: Global smoothness.}

We will now show that $\phi$ is $C^{\infty}$. The manifold $M$ is
foliated by the  smooth foliations
${\cal W} ^- _H$ and ${\cal O}$ whose tangent bundles are transverse
and span the $TM$. We showed above that $\phi$ is $C^{\infty}$ along
all the foliations  ${\cal W} ^- _H$. As $\phi$ is already
$C^{\infty}$ along ${\cal O}$ by construction of $\phi$, then $\phi$
is smooth along a full set of directions. Now it follows from standard
elliptic operator theory that $\phi$ is a
$C^{\infty}$-function (as the elliptic operator $\frac{\partial
^{2s}}{\partial x_1 ^{2s}} +
\ldots  \frac{\partial ^{2s}}{\partial x_n ^{2s}}$ sends $\phi$ to a
continuous function where
$x_1,\ldots, x_n$ are coordinates subordinate to the foliations above.)
Global $C^1$ regularity in the finite differentiability case 
directly follows from the $C^1$ regularity along each foliation.

This finishes the proof of Theorem~\ref{thm-main}. \QED

\subsection{Proof of Corollaries to Theorem 1}

\subsubsection{Proof of Corollary 2}

It suffices to check the conditions of Theorem~\ref{thm-main}. 
Since the maximal nontrivial intersections of stable manifolds are exactly
coarse Lyapunov spaces and all one-parameter subgroups are ergodic, we can 
 pick $t \, b$ to be any  one-parameter subgroup  in the boundary of the
 Lyapunov halfspace.

\subsubsection{Proof of Corollary 3}
The argument for the Weyl chamber flows has been given 
 before the statement of the corollary in Section 2.1.

 We will use the notation established
in the description of twisted Weyl chamber flows before Corollary 3. It suffices to
prove that all one-parameter subgroups of $A$ are ergodic. If $G$ is simple, 
 as above, Moore's theorem yields the ergodicity of all
one-parametere subgroups so that we can apply  Corollary 2 above. For the
general case, we will need to use the results of Brezin and Moore.
 Let us briefly
recall some notions from  \cite{Brezin-Moore}. Call a  simply connected Lie group
{\em Euclidean} if it is the universal cover of an extension of a
vector group by a compact abelian group. A solvmanifold is called {\em
Euclidean} if it is the quotient of a Euclidean group by a closed subgroup.
Further, given a connected Lie group $H$ and a closed subgroup $D$, we call the
pair $(H,D)$ {\em admissible} if $D$ normalizes a closed solvable subgroup
containing the radical of $H$. Note that the pair $( G\ltimes  N, \Gamma
\ltimes  \Lambda )$ is always admissible.
Brezin and Moore show that every admissible pair $(H,D)$ has a maximal
Euclidean quotient $H/L$ where $D \subset L$. Let $H^{(n)} =[ H^{(n-1)} ,
H^{(n-1)}] $ denote the $n$-th element of the derived series where $H^{(0)} =H$
. Then $ H^{(n)}$ will eventually be constant.
Let $H^{(\infty)}$ denote this subgroup. Brezin and Moore show that $
H^{(\infty)} \subset L$. For our special case of $H = G \ltimes  N$, where
the twisting homomorphism contains an Anosov element, the derived series is
constant right away, as follows from  a consideration of eigenvalues. Hence the
maximal Euclidean quotient is a point.

Brezin and Moore show that a one parameter subgroup of $H$ acts ergodically on
$H/ D$ for an admissible pair $(H,D)$ if and only if the quotient flows on the
maximal Euclidean quotient and the maximal semisimple quotient are ergodic.The
latter quotient is obtained by factoring $H$ by the closure of $D \, R$ where
$R$ is the radical of $H$. In our case, the latter is isomorphic to $G/\Gamma$
as $ \Gamma\ltimes  N$ is closed in $G\ltimes  N$.
Since $\Gamma$ is irreducible, all one- parameter subgroups of $A$ are ergodic
on $G/\Gamma$. Since the maximal Euclidean quotient is a point, we get
ergodicity of one-parameter subgroups of $A$.

{\it Remark.} A.Starkov pointed out to us  an alternative argument, avoiding the 
discussion of admissibility, which uses Theorem 8.24  from \cite{Ragh}  and 
certain results proved independently by himself and by D.Witte.

\subsubsection{Proof of Corollary 4}

Suppose $\Zk$ acts on a compact nilmanifold $M$.  Then the suspension
$(M \times \Rk) /
\Zk$ factors over $T^k :=\Rk / \Zk$ with fiber $M$. Since the strong
stable manifolds in the suspension of any $b \in \Rk$  are contained
in the fibers, the condition on intersections of unstable manifolds in
Theorem~\ref{thm-main} is satisfied if we know that the ergodic
components of any one-parameter subgroup
$t\, a$ of the suspension consist of whole fibers. It clearly suffices
to check this on some finite cover of $M$. Thus, we will assume that
$M=N/ \Lambda$ is a nilmanifold, and $N$ is a simply connected
nilpotent group. 

Note that the linearization of the $\Zk$-action extends (in general, not uniquely) to
$\Rk$. Let $S= \Rk
\ltimes  N$ be the solvable group defined by this extension. Then 
$\Sigma := \Zk \ltimes  \Lambda   $ is a lattice in $S$, and the
natural fibration  of $S/ \Sigma $ over $T^k$ is just the fibration
defined by the suspension. Furthermore, the action of a one-parameter
subgroup in the suspension is just  the homogeneous action of a
one-parameter subgroup of $\Rk$ on $S/ \Sigma $ via the embedding of
$\Rk$ into $S$. 

As in the proof of the previous corollary we will now  apply the results of
 J. Brezin and
C. C. Moore on homogeneous flows \cite{Brezin-Moore}. For a general compact
 solvmanifold, Brezin and
Moore show that there always exists a maximal quotient which is
Euclidean. Moreover, this quotient is unique and is called the  {\em
maximal Euclidean quotient}. They further show that a homogeneous flow
on a solvmanifold is ergodic if the quotient flow on the maximal
Euclidean quotient is ergodic \cite[Theorem 6.1]{Brezin-Moore}.

For the case of our suspension flow, our claim about ergodic
components will follow quickly once we determine the maximal Euclidean
quotient of suspensions.

\begin{Lemma} The maximal Euclidean quotient of the suspension
$S/\Sigma$ of an algebraic 
$\Zk$-action for which every nontrivial element is weakly mixing is
$T^k$.
\end{Lemma}

Assume the lemma. Let $T^r$ be an ergodic component of the quotient
action on
$T^k$ of some one-parameter group $ t\, a \in \Rk$. Then $T^r = \R^r
/\Z ^r$ for $\Z^r \subset
\Zk$. Then the preimage of $T^r$ in $S/\Sigma$ is a solvmanifold
$S^r/\Sigma ^r$, and in fact is the suspension of the $\Z ^r$ action
on $M$.  Applying the lemma to $S^r/\Sigma ^r$, the flow on
$S^r/\Sigma ^r$ is ergodic by Brezin's and Moore's theorem. Hence the
ergodic components of the one-parameter group on $S/ \Sigma$ are all
preimages of subtori of $T^k$, and in particular contain whole fibers.

{\em Proof of the Lemma}\,:  Let $S/D$ be the maximal Euclidean
quotient of $S/
\Sigma$. Since $T^k$ is a Euclidean quotient,  $D \subset 
\Zk \ltimes  N$. We are done if  the connected component of the
identity $D^0$  of $D$ is  $N$. Suppose not. Then $D^0$ does not
project onto $N/ [N,N]$  since a subspace of the Lie algebra of $N$
complementary to the commutator algebra generates the whole Lie
algebra. Also note that the linear part of the 
$\Zk$-action on $N/D^0$ has eigenvalues of modulus 1 since they
correspond to eigenvalues of the flow on the maximal Euclidean
quotient. 

Since $D \cap N \supset \Sigma \cap N = \Lambda$, $D \cap N$ is
cocompact, and hence Zariski dense in $N$ 
(cf. e.g. 
\cite[Theorem 2.3]{Ragh}).  Thus $D^0 \subset D \cap N$ is normal in
$N$. Since $\Zk \subset D$, $\Zk$  normalizes $D^0$. Since the  
normalizer
 of 
$D^0$ is algebraic it contains a proper $R^k$ in which our $Z^k$ is a lattice.
We conclude that $D^0$ is normal in $S$.  Hence $D/D^0$ is a lattice
in $S/D^0$, or equivalently, $\Zk \ltimes  (D \cap N) /D^0$ is a
lattice in  $\Rk \ltimes  N/ D^0$. This implies that 
 $\Z ^k$ leaves the lattice $(D \cap N)/D^0$ in $N/D^0$ invariant. 
 It follows that $\Z ^k$ leaves the lattice $ (D \cap N)/ D^0 \,
[N,N]$ in 
 $N/ D^0[N,N]$ invariant. Since $\Zk$
 acts on $N/D^0[N,N]$ with linearization with eigenvalues of modulus 1,
these
 eigenvalues have to be roots of unity. Thus the linearization of
$\Zk$ on $N/[N,N]$ also has roots of unity as eigenvalues unless
$D^0=N$. By Parry's characterization of weakly
 mixing automorphisms of nilmanifolds, no nontrivial element of $\Zk$
can have
 eigenvalues of modulus 1 on  $N/[N,N]$, which is a contradiction
\cite{Parry}.

\subsubsection{Proof of Corollary 5}

 First we establish that any $C^1$-small perturbation of  a standard Anosov action is
$C^{\infty}$-orbit equivalent to the original action. In the case of the suspension of an
action by automorphisms of an infranilmanifold this follows from Corollary 4 and in the
cases of Weyl chamber flows and twisted Weyl chamber flows  from Corollary 3. The
remaining extension cases are handled similarly and we omit the necessary but somewhat 
tedious algebraic arguments. Hence the perturbed action is smoothly conjugate to a
$C^{\infty}$-time change of the standard Anosov action. However, we showed in \cite{KS1}
that such time changes are 
$C^{\infty}$-conjugate to the original action up to an automorphism. 

Furthermore, as the orbit equivalence in question can be
chosen $C^1$-close to the identity, the resulting time change is close
to the identity.  One can further  choose the conjugacy between the
time change and the original action close to the identity. In fact,
we  solve a cocycle equation in \cite{KS1} by defining the coboundary
$P$ as a sum of the cocycle over the forward orbit. It follows  easily
that the derivative is small along stable manifolds, and similarly
along unstable manifolds. Hence the coboundary is close to $0$ in the
$C^1$-topology if it is close to $0$ at some point. As we can always
pick the coboundary to have 0 average with respect to the invariant
volume,
$P$ is close to 0 somewhere. That the automorphism $\rho$ is close to
the identity and unique follows as it is an average of the cocycle
determined by  the time change  over the given volume. 

Finally, suppose that $\phi _1$ and $\phi _2$ are two conjugacies
close to the identity.  Then 
$\phi  =\phi _2 ^{-1} \, \phi _1$ commutes with the $\Rk$-action and
is close to the identity. Since  
the identity is an orbit
equivalence of the $\Rk$-action with itself, $\phi$ has to take orbits
to themselves by the uniqueness part of Theorem~\ref{thm-main}. Now
fix a point $p$ with a dense
$\Rk$-orbit. This is possible as the action is ergodic.  Then find the
translation $a \in \Rk$ such that $a\, p = \phi (x)$. Then $a$ and
$\phi $ coincide on the $\Rk$-orbit of $p$, and hence everywhere as
that orbit is dense.
 \QED

\subsection{The reductive case}

Let us make a few remarks about the changes needed to prove Theorem 6
and Corollary 7. Suppose we have an Anosov action of a reductive group
$H = L\, C$    on $G/\Lambda$ where $G$ is a semisimple Lie group, $C$
is a subgroup of a split Cartan, and $\Lambda \subset G$ is a lattice.
First note that the orbit foliation ${\cal O}$ of $H$ is the neutral
foliation of a suitable element $a \in C$. Let ${\cal R}$ be a small
perturbation of ${\cal O}$. Perturbing vectorfields from $C$ as
before, we can recover  ${\cal R}$ as a normally hyperbolic invariant
foliation for some  perturbed vector field. By structural stability,
we get an orbit equivalence. Now we proceed as in the abelian case,
and define a continuous action of $C$ which is subordinate to ${\cal
R}$. The difference to the abelian case is that this action is not
transitive on  ${\cal R}$ anymore. Still, it is  ergodic and hence
topologically transitive. Hence we get  transitive groups on the
leaves of the  coarse Lyapunov foliations, as before. The rest of the
arguments apply verbatim.

Corollary 7 follows the same way, as in fact the group $C$ still acts
on 
$\Gamma \setminus G / M_P$.

\subsection{Automorphisms of Weyl chamber flows}

\def \Ad
 {\mbox
{ Ad } } For an action $\alpha$ of a group $G$  on a manifold $M$ let the 
 {\em automorphism group}  $Aut (\alpha)$ be the centralizer of 
$\alpha (G)$ in the homeomorphism group of $M$. 

In this subsection  we prove  a result showing that a Weyl chamber
flow is a maximal abelian action on its ambient manifold.

\begin{Theorem}   If  $\alpha$ is a   Weyl chamber flow on
$K\setminus H / \Lambda$ where $H$ is a semisimple Lie group  of rank 
$k\geq 2$ without
compact factors, then the   automorphism group of
$\alpha$ consists of $C^{\infty}$-diffeomorphisms and is a finite
extension of $\Rk$, i.e. the $\Rk$ 

action itself is a 
subgroup of finite index in the centralizer. 
\end{Theorem}

This result has more to do with semisimplicity of $H$ than with the
rank assumption which is central for all the rest of our results. For
example, the assertion  is also well-known for the Anosov actions of
the split  Cartan subgroups of rank one semisimple Lie groups $G$ of
the non-compact type on a compact quotient $M
\setminus G / \Gamma$ where $\Gamma$ is a cocompact lattice and $M$ is
the compact part of the centralizer of a split Cartan subgroup.
        
On the other hand, if $\alpha$ is an Anosov action on a nilmanifold,
then all its automorphism are diffeomorphisms as is well known but
the centralizer naturally may be much bigger than the action.      
 The situation is similar for twisted Weyl chamber
flows which are toral fiber bundles over a Weyl chamber flow.  We will
not discuss  this case here however.

\PROOF Let   $ \alpha$ be a Weyl chamber flow. 
 Let  $ N = K \setminus H /\Lambda $ be the space  $\alpha$ acts on.  
Suppose that $\phi \in Aut(\alpha)$. Let $\Gamma = \pi _1 N$, the
fundamental group of $N$.  Then
$\phi$ lifts to a 
$\Gamma$-equivariant map on the universal cover  $\tilde{N}$ of $N$. 
Note that $\tilde{N} $ is a homogeneous space of the universal cover
$\tilde{H}$ of $H$ and let $\rho: \tilde{H} \rightarrow \tilde{N}$ be
the   factor map. Lift $\Rk$ to $\tilde{ H}$. The space of weak stable
manifolds of 
$a \in \Rk$ on $\tilde{N}$  is isomorphic to the space of weak stable
manifolds of $a \in \Rk$ on
$\tilde{H}$. The latter can be described as the set of parabolics $H/P
^a$ where the  Levi component  of $P ^a$ is  the centralizer
 of $\log a$ in $\tilde{H}$. Thus we can map $x \in \tilde{N}$ to 
the  parabolic $\sigma _a (x)$ which is the stabilizer of the weak
stable manifold 
$\rho ^{-1} (W^{ws} _a (x))$ of $a$ at  any point in $\rho ^{-1} (x)$. 

Note that  for all $g \in \tilde{G}$, $\sigma _a (x\,g) = g \,\sigma
_a (x)\,
 g^{-1}$. Hence $\sigma _a$ is surjective onto $H/P ^a$.  If $\sigma
_a (x) = 
\sigma _a (y)$ then $y \in W^{ws} _a (x)$ as follows from the
equivariance above and the self-normalizing property of $P ^a$. 
Therefore and because $\phi$ takes weak stable manifolds of
$a$ to weak stable
 manifolds of $a$ we can define a map $\bar{\phi} _a: H/P ^a
\rightarrow H/P ^a$ by
\[\bar{ \phi} _a \circ \sigma _a = \sigma _a \circ \tilde{\phi}.\]
 
 Let $a$ belong to an open face $F$ of a Weyl chamber. If $b$ belongs
to the
 boundary of $F$ then 
\[ \sigma _a (x) \subset \sigma _b (x).\] Thus if $P_1 \subset P_2$
are parabolics in $H/ P_a$ and
$H/ P_b$ respectively, then 
\[ \phi _a (x) \subset \phi _b (x).\] Therefore the $\phi _a$ define
an automorphism of the Tits building  $\Delta$  attached to
$\tilde{H}$ where $a$ ranges over a closed Weyl chamber.

We will show now that  $\sigma _a$  and  $\phi _a$ are independent of
$a$ for all regular elements
$a \in \Rk$. First note that they do not change as
$a$ varies over an open  face of a Weyl chamber since the weak stable
manifolds stay the same. If
$a$ and  $a'$ belong to  adjacent open Weyl chambers ${\cal C}$ and
${\cal C'}$ then $\phi _a = \phi _{a'}$.  Indeed, 
$\phi _a$ and  $\phi _{a'}$  coincide on the faces of the type of the
intersection of $\bar{{\cal C}}$ and $\bar{{\cal C'}}$ by
equivariance. Thus they coincide on all of $\Delta$ 
 by a convexity argument.  Since 
 any two Weyl chambers are connected by a chain of adjacent Weyl
chambers,  it follows that all
$\phi _a$ coincide.

Thus $\tilde{\phi}$ defines a unique automorphism  $\bar{\phi}$ of
$\Delta$  which is clearly continuous for the topology on $\Delta$
inherited from 
$\tilde{H}$. Note that $\tilde{\phi}$ is completely determined by
$\bar{\phi}$ transversally to the
$\Rk$-orbits by projecting $x \in \tilde{N}$ simultaneously to $\sigma
_a (x)$ and $\sigma _{-a} (x)$ for a regular element $a$. Since a  
topological automorphism of $\Delta$ is $C^{\infty}$ (for the smooth
structure
 on $\Delta$ induced by the $\tilde{G}$-action) \cite{BS0} and the
maps 
$\sigma _a$ and $\sigma _{-a}$ are $C^{\infty}$ we see that
$\tilde{\phi}$ is
 $C^{\infty}$ transversally to the $\Rk$-orbits. As $\tilde{\phi}$ is 
$C^{\infty}$ along the orbits, $\tilde{\phi}$ is $C^{\infty}$. 

 Consider the group $\tilde{Aut} (\alpha)$ of lifts of automorphisms of
 $\alpha$ to $\tilde{N}$ and the map $\tilde{Aut} (\alpha) \rightarrow
 Aut (\Delta)$ where $Aut (\Delta)$ is the group of topological
automorphisms
 of $\Delta$. Then $\Rk$ embeds canonically into $\tilde{Aut}
(\alpha)$, and
 is the kernel of the map into $Aut (\Delta)$. Indeed, if $\bar{\phi}
= id$  then   $\tilde{\phi}$ maps every $\Rk$-orbit in $\tilde{N}$
into itself as 
$\tilde{\phi}$ is transversally determined by $\bar{\phi}$. Then we
can write
 uniquely $\tilde{\phi} (x) = a(x) \, x$ for $ x \in \tilde{N}$, and
thus 
$\phi (x) = a(x) \,x$ at least for all $x \in N$ without isotropy in
$\Rk$.
 Since  $\phi$ is an automorphism, $a(x)$ is constant along  orbits,
and  hence constant by ergodicity of the $\Rk$-action.   By 
\cite{BS0}, the topological automorphism group of a topological 
building is a finite extension of $\Ad G$. Hence the image of
$\tilde{Aut}  (\alpha)$ in $Aut (\Delta)$ contains the intersection 
$I$ with $\Ad H$ as
 a subgroup of finite index. If $\bar{\phi} \in I$ then $\bar{\phi}$
belongs  to the normalizer of
$\Ad \Lambda$ in $\Ad H$ since $\bar{\phi}$ is $\Lambda$  equivariant. 
Further, if $\bar{\phi} \in
\Ad \Lambda$ then $\phi = id$. As  the group of outer automorphisms of
$\Lambda$ is finite (as follows e.g. from Mostow's rigidity theorem),
the theorem is proved. 

\section{Algebraic Anosov  actions of   lattices}

We will use our results on smooth local rigidity of  affine 
Anosov actions of abelian groups to investigate  algebraic Anosov 
actions of lattices in semisimple Lie groups of higher rank. In particular,
we obtain local rigidity for general affine Anosov actions of such lattices
in  Theorem~\ref{thm-lattice-anosov} below. 

We will first prove a general result, Proposition~\ref{prop-locrig}, asserting 
the local  rigidity of certain 
affine lattice actions modulo the local rigidity of  suitable large Abelian
subgroups.   Theorem~\ref{thm-lattice-anosov} itself  follows 
albeit not so quickly from this and our results in Section 2. The
main two steps  involve verifying the assumptions of 
Proposition~\ref{prop-locrig}  and passing from rigidity on a
subgroup of finite index to the rigidity of the whole action.
Proposition~\ref{prop-locrig} may  prove  useful in establishing 
local rigidity results beyond the Anosov case.

Let $G$ be  a linear semisimple  Lie 
group without compact factors all of whose simple factors have real rank at 
least 2. Let $\Gamma$ be an irreducible lattice in $G$.
Call a subgroup $\Delta \subset \Gamma$ a {\em Cartan 
subgroup}
of $\Gamma$ if there is a Cartan subgroup $D$ in $G$ 
with maximal $\R$-split component such that $\Gamma$ intersects $D$ in a
lattice in $D$. In particular, $\Delta$ is an abelian group whose rank is
the real rank of $G$. By a theorem of G. Prasad and M. Raghunathan, any lattice
$\Gamma$ contains Cartan subgroups \cite[Theorem 2.8]{Prasad-Ragh}.

Let $M$ be a manifold of dimension $N$ with a \ci -action $\rho$ of a group 
$\Gamma$ preserving a probability measure $\mu$.
We call a measurable framing $\tau$ of a manifold $M$ a {\em superrigidity 
framing} for $\rho$ if there is a homomorphism $\sigma : \Gamma \rightarrow
GL(N,\R)$, a compact subgroup $C \subset GL(N,\R)$ commuting with the image of
$\sigma$ and a measurable cocycle $\kappa: \Gamma \times M \rightarrow C$ such
that for all $\gamma \in \Gamma$ and $\mu$- a.e. $x \in M$
\[ D \rho (\gamma) (\tau (x)) = \tau (\gamma \, x) \, \sigma (\gamma)
\, \kappa (\gamma, x) .\] 
Then $\pi$ is called a {\em superrigidity representation}, and $\kappa$ a 
{\em superrigidity cocycle}. Furthermore, we  call a framing $\tau$ on a 
homogeneous space $L/
\Lambda$ {\em translation invariant} if $\tau$ lifts to  a framing of $L$
invariant under right translation by $L$.

\begin{Proposition}      \label{prop-locrig}
Let $G$ be  a linear semisimple  Lie 
group without compact factors all of whose simple factors have real rank at 
least 2. Let $\Gamma$ be an irreducible lattice in $G$,  
$\Delta$  a Cartan subgroup of $\Gamma$, and  
$\rho$ an ergodic affine algebraic action of $\Gamma$ on a manifold $M$. Assume 
that the restriction of $\rho$ to $\Delta$ is locally \ci -rigid, and that we
can choose the conjugacy $C^1$-close to the identity. Assume further that 
for some $a \in \Delta$, $\rho (a)$ as well as all lifts of finite powers
 to finite connected covers 
are weakly mixing.  Then any \ci -action $\tilde{\rho}$ of $\Gamma$,
sufficiently close to
$\rho$ in the $C^1$ -topology, is \ci -conjugate to an action with a
translation invariant superrigidity framing for a subgroup $\Gamma _0$ of
$\Gamma $ of finite index.
\end{Proposition}

Let us emphasize that the affine lattice action need not be Anosov in this
proposition.  
\vspace{.5em}

\PROOF 
The basic idea of the proof is to compare a superrigidity framing for
the perturbed action with a translation invariant framing using the
dynamics  of the Cartan subgroup $\Delta$ and  see how the framings get
transformed under 
$\Delta$. In effect, we prove a certain uniqueness of superrigidity 
framings for suitable affine actions of $\Delta$.

\vspace{1em}

{\bf Step 1: Measurable superrigidity} 
  We  first need to review some aspects of R. Zimmer's   superrigidity
theorem  for cocycles \cite{Zimmer}. Suppose $\tilde{\rho}$ is a 
$\Gamma$-action on  a manifold 
$M$ preserving a finite measure on $M$.  As $G$ is linear, 
$G$ is a quotient of the maximal algebraic factor  $\bar{G}$ of the 
universal cover of $G$. Lifting $\Gamma$ to a lattice in
$\bar{G}$, we may assume without loss of generality that $G$ is the 
maximal algebraic factor of its universal cover. Let $P \rightarrow M$
be a principal bundle over $M$ with structure group $L$ on which
$\Gamma$ acts  smoothly by bundle automorphisms covering the action on
the base. We will denote this extension $D\tilde{\rho}$ since in our 
application it will be generated by the differentials of a smooth
action in the base.  The {\em algebraic hull}  of
$P$ and
$\tilde{\rho}$ is the smallest algebraic subgroup $H$ of $L$ (up to
conjugacy), for  which there is a measurable section $\tau$ of $P$
which transforms under
$\tilde{\rho}$ by some elements $\alpha (\gamma,x)$ of $H$: \[D
\tilde{\rho} (\gamma) (\tau (x)) = 
\tau ( \tilde{\rho} (\gamma) (x))\: \alpha (\gamma,x).\]
 Then $\alpha$ is a cocycle over $\tilde{\rho}$ with values in $H$. 
Furthermore, $H$ is always reductive with compact center as was shown
for uniform lattices by R. Zimmer in \cite{Zimmer1} and 
for  non-uniform lattices by J. Lewis in \cite{Lewis}. In order to be
able to apply  Zimmer's superrigidity theorem  for cocycles
\cite{Zimmer}, the algebraic hull needs to be (Zariski)  connected.
This can be achieved by passing to a suitable finite cover  
$M'$ of $M$ as follows. Let $H^0$ be the Zariski connected component
of $H$. Let $P$  be the reduction of the frame bundle of $M$ with
structure group $H$. Set 
$M' =\frac{P \times H^0 \setminus H}{H}$. Then $M' \rightarrow M$ is a 
finite  cover,  on which $\Gamma$ acts as a quotient of the bundle 
automorphisms on $P$ by $H$. Note that this lift of the
$\Gamma$-action  preserves a finite volume. Moreover, $P \rightarrow
M'$ is a principal 
$H^0$-bundle on which $\Gamma$ acts via bundle automorphisms. 
Measurably, $M'$ is $M \times H/H^0$. The action is a skew product
action.  Let $\beta$ be the projection of $\alpha$ to $H/H^0$, and
define the skew  product action of $\Gamma$ on $M'$ by 
\[ \gamma (x,[h]) =(\gamma \, x, \beta (\gamma ,x) [h]).\] Note also
that the skew product leaves the product of the invariant  volume on
$T^N$ and counting measure invariant.  Now the following theorem is
the desired corollary of measurable superrigidity (cf. e.g.
\cite[Theorem 3.1]{Qian-Zimmer}).

\begin{Theorem}    \label{thm-superrigid}
Passing to the finite cover $M'$ of $M$ described above, there is a 
measurable section
 $\tau ^*$ of $P$ which transforms under $\tilde{\rho}$ by a cocycle of 
the form $\pi
 (\gamma)  \: \kappa (\gamma, x)$ where $\pi: G \rightarrow H^0$ is a 
 homomorphism, and 
 $\kappa$ is a cocycle taking values in a compact normal subgroup of 
$H^0$
 centralizing $\pi (\Gamma)$.
  \end{Theorem}
\vspace{1em}

Let us note here that  $M'$ is only constructed in the category of
measure spaces since the reduction $P$ of the frame bundle with
structure group $H$  a priori is only measurable. This will cause some
complications in the next step of the proof.

\vspace{1em}

{\bf Step 2: Regularity of the superrigidity framing}
We assume the notations from the statement of 
Proposition~\ref{prop-locrig}.  As $\rho$ is affine,
$M$ is the quotient of a connected Lie group $L$, $M = L / \Lambda$. 
We will show that if
$\tilde{\rho}$ is sufficiently close to $\rho$ in the $C^1$-topology,
then 
$\tilde{\rho}$ is \ci -conjugate to $\rho$. How close $\tilde{\rho}$
needs  to be to $\rho$ will be specified during the proof. Since the
restriction
 of $\rho$ to $\Delta$ is locally $C^1$-rigid, and since  we may
choose the
 conjugacy $C^1$-close to the identity, we may assume that
 $\rho$ and $\tilde{\rho}$ coincide on the Cartan subgroup $\Delta$.  

Since $\Gamma$ has Kazhdan's property and $\rho$ preserves Haar
measure, a perturbation 
$\tilde{\rho}$ sufficiently close to $\rho$ in the $C^1$-topology
preserves an absolutely continuous probability measure (cf. e.g.
\cite[Lemma 2.6]{KLZ} or \cite{Seydoux}).

Now let us apply Theorem~\ref{thm-superrigid} to $\tilde{\rho}$ and
its  extension by derivatives to the frame bundle of $M$. In
particular, we  let $M'$ denote the finite cover of $M$ for the
$\tilde{\rho}$ action
 and $\tau ^*$ be a superrigidity framing  of $M'$ as in
Theorem~\ref{thm-superrigid}. 
Denote the action 
of $\Gamma$ on $ M'$ extending $\tilde{\rho}$ by $\hat{\rho}$. 

Let $\tau $ be a translation invariant framing of $M$.  
We will now show that we can pick a superrigidity framing  
 for a subgroup of finite
index of $\Gamma$ on $M'$ which is just a
constant translate of $\tau$ modulo a compact subgroup.

Since $M'$ a priori is only a measurable gadget and our arguments
involve some topology, we will need to argue  on the base manifold
$M$. Thus we will consider the given superrigidity  framing $\tau ^*$
as a finite collection of frames $\tau ^* (x)$ over any point $x$ in
$M$ of cardinality $c=\# (H/H^0 )$ corresponding to the $c$ points in
$M'$ which  cover $x$. Let ${\cal G}$ denote the space of unordered
$c$-tuples of points in
$GL(N,\R)$. Then  we can  define a measurable   function 
$\beta: M \rightarrow {\cal G}$ such that $\tau (x) \, \beta(x) = 
\tau ^* (x)$. We can think of $\beta$ as a multivalued function with
values in $GL(N,\R)$.

  Since the restriction of $\hat{\rho}$  to $\Delta$ is just the 
original action $\rho$ on $M$, the  framing $\tau $ is transformed
under $\Delta$ by $\sigma$.  To  
simplify   
notations, we will denote the
action of $\hat{\rho} (b)$  by $b\, x$  for $b \in \Delta$. By
assumption, there is an element $a \in \Delta$ which is weakly mixing
on $M$. We will also use the symbol $d$ for the bundle extension.

Since the translation invariant framing $\tau$ is transformed  under
$\Delta$ by the automorphism  $\sigma$ , we get for
$a \in \Delta$ and $x \in M$:
\[ d\, a (\tau ^* (x)) = d\, a (\tau (x) \, \beta (x)) =\tau (a\, x)
\: 
\sigma (a) \: \beta (x) \] where  $g\in GL(N,\R)$  acts on ${\cal G}$
in two ways by either multiplying each element of the $c$-tuple by $g$
on the left or the right. 
 On the other hand, by Theorem~\ref{thm-superrigid} we have for a.e.
$x $,
\[ d\, a (\tau ^* (x)) = \tau ^* (a\, x)\: \pi (a)\: \kappa (a, x)=
\tau (a\, x)\: \beta (ax) \: \pi (a)\: \kappa (a, x) .\] Hence we get
for a.e. $x$ and all $a \in \Delta$,
\[ \sigma (a) \, \beta (x) = \beta (ax) \, \pi (a)\, \kappa (a, x) .\]
By Lusin's theorem, for every $\varepsilon >0$, there exists a closed 
subset
$M_{\varepsilon} \subset M$ such that the measure of $M_{\varepsilon}$ 
is at least $1- \varepsilon$ and $\beta$ is uniformly continuous on
$M_{\varepsilon}$. Pick $a \neq 1 \in \Delta$. As    the
transformation induced by $a$ on  $M_{\varepsilon}$ preserves the
induced measure, a.e. point  $x \in M_{\varepsilon}$ is recurrent.
Consider such a point $x$. Then there is a sequence of integers $n_k
 \rightarrow \infty$ such that $a^{n_k} \, x \in M_{\varepsilon}$ and 
 $a^{n_k} \, x \rightarrow x$. Let $b _0 $ be an element of $\beta
 (x)$. Then there are elements $b_{n_k} \in \beta (a^{n_k} (x))$ such
that
 $\sigma (a) ^{n_k} \, b_0 = b_{n_k} \, \pi (a)^{n_k} \, 
 \kappa (a^{n_k} , x)$. 
 Since the values of $\kappa$ and the image of
 $\pi$ commute, it follows that
 \[b_0 ^{-1}\, \sigma (a ^{n_k}) \,  b_0 = b_0 ^{-1} \, 
 b_{n_k} \, \kappa (a ^{n_k}, x) \, \pi (a ^{n_k}) .\]
 Since $\beta (a^{n_k} \,x)$ converges to $\beta (x)$, it follows that 
 $b_0 ^{-1} \,b_{n_k} $ form a bounded sequence of matrices. 
 Moreover  $\kappa (a ^{n_k}, x)$ lie in a compact group.
 Therefore,  the eigenvalues of $\sigma (a)$ and 
 $\pi (a)$ have the same absolute values. In other words, they have
the same
 Lyapunov exponents. Moreover the dimensions of the sum
 of eigenspaces with fixed absolute value (the Lyapunov spaces)
coincide.
 Since $a$ is semisimple, both $\sigma (a)$ and $\pi (a)$ are
semisimple.
 Hence we can write them as commuting products 
 $\sigma (a) = \sigma _c (a) \, \sigma _{nc} (a)$ and 
 $\pi (a) = \pi _c (a) \, \pi _{nc} (a)$
 where $\sigma _{nc} (a)$ ($\pi _{nc}(a)$) has only positive
eigenvalues and 
 $\sigma _c (a)$ ($\pi  _c (a)$) has only eigenvalues of modulus 1. It 
 follows from the discussion above that $\sigma _{nc} (a)$ and $\pi
_{nc}(a)$
 are conjugate in $GL(N,\R)$. Let us also note that these
decompositions are unique.
 Hence, any element commuting with $\sigma (a)$ or $\pi (a)$ also
commutes
 with their components.
 
 Pick an element $e \in GL(N,\R)$ such that 
$e \, \sigma _{nc} (a) \, e^{-1} = \pi  _{nc}(a) $. Set 
$A= e^{-1}\, \pi _c (a) \, e$ and  $\bar{\kappa} (a,x) \df e^{-1} \,
\kappa (a,x) \, e$.  Note that  $\sigma _{nc} (a) = e^{-1}\, \pi _{nc}
(a) \, e$ commutes with both
$A$ and the images of $\bar{\kappa}$. Also set $\bar{\tau} = 
\tau ^*  \,e$. Then $\bar{\tau}$ is a measurable superrigidity framing 
on $M'$ for $\hat{\rho}$ for all of $\Gamma$ for the representation
$e^{-1} \, \pi \, e$ and the  cocycle $\bar{\kappa}$. Note that
$\bar{\tau}$ gives rise to a measurable reduction of the frame bundle
which is just a constant translate of $P$. Hence we may assume that
$P$ is the measurable reduction determined by $\bar{\tau}$. Again, we
will interpret $\bar{\tau} (x)$ as an unordered $c$-tuple of frames
for $x$ in the base manifold $M$.  Let $\psi : M
\rightarrow {\cal G}$ be the measurable function such that 
$\bar{\tau} (x) =\tau (x) \psi (x)$ for $x \in M$.  Then we get for
all $n \in \Z$
 \[ d\, a^n (\bar{\tau} (x)) =d\, a^n ( \tau ^*( x) \, e )
	=\tau ^*(a^n \,  x)\, \pi (a^n) \kappa (a^n,x) \, e
	= \bar{\tau} (a^n\, x)\, e^{-1} \, \pi (a^n)  \, \kappa (a^n,x) \,
	e\] 
\vspace{-1.5em}
	\[=\bar{\tau} (a^n\, x)\, e^{-1} \, \pi _{nc} (a^n) \, e\, 
	e^{-1} \pi _c (a^n) \, e \, \bar{\kappa} (a^n,x)
	= \tau (a^n \, x) \, \psi (a^n \, x) \sigma _{nc} (a^n) \, A^n \, 
	\bar{\kappa} (a^n,x).\]
	
We also get 
	\[d\, a^n (\bar{\tau} (x)) =  d\, a^n (\tau (x) \,\psi (x)) =
\tau (a^n \, x) \, \sigma (a^n) \psi (x).\] Therefore  we see that
\[ \psi (a^n \, x) \, \sigma _{nc} (a) ^n \, A^n \, \bar{\kappa}
(a^n,x) =
\sigma _c (a) ^n \, \sigma _{nc} (a) ^n \, \psi (x).\] Let $B= \sigma
_c (a)$. Since $\sigma _{nc} (a)$ commutes with $A$  and
$\bar{\kappa}$, we get
\[B^{-n} \, \psi (a^n \, x) \, A^n  \, \bar{\kappa} (a^n,x) = 
\sigma _{nc} (a) ^n \, \psi (x) \, \sigma _{nc} (a) ^{-n}.\]

By Lusin's theorem, for every $\varepsilon >0$, there exists a closed 
subset
$M_{\varepsilon} \subset M$ such that the measure of 
$M_{\varepsilon}$ is at least $1- \varepsilon$ and $\psi$ is 
continuous on
$M_{\varepsilon}$.   The transformation induced by $a$ on 
 $M_{\varepsilon}$ is measure preserving. Hence, by Poincar\'{e}'s 
 recurrence theorem,  for a.e. $x \in M_{\varepsilon} $ there are 
sequences of integers $n_k
\rightarrow \infty$ and $m_k \rightarrow -\infty$ such that 
$a^{n_k} \, x, a^{m_k} \, x \in M_{\varepsilon}$, $a^{m_k} \, x 
\rightarrow x$  and  $a^{n_k} \, x \rightarrow x$ as $k \rightarrow
\infty$. Since $A$, $B$ and the range of $\bar{\kappa}$ belong to
compact groups, it follows  that
$\sigma _{nc} (a) ^{m_k} \,\psi (x)\, \sigma (a)  _{nc} ^{-m_k}$ and
also 
$\sigma (a) _{nc} ^{n_k} \,\psi (x)\, \sigma  _{nc}(a) ^{-n_k}$ stay
bounded as $k 
\rightarrow \infty$.  As $\sigma  _{nc}(a) \in GL(N,\R)$ is
diagonalizable, this is only possible if 
$\psi (x)$ commutes with $\sigma  _{nc}(a)$. Thus we get that
\[\hspace{15em} \psi (a\, x) \, A \, \bar{\kappa} (a,x) = B \, \psi (x)
 \hspace{15em} (*)\] for a.e. $x \in M$ as $M_{\varepsilon} $ has
measure arbitrarily close to 1.

Now we pick the element $a$ satisfying the weak mixing assumption and
will  show that this last cohomology equation forces $\psi$ to be
constant modulo a suitable compact subgroup. Let ${\cal A}$ and ${\cal
B}$ denote the closure of the powers
$A^n, n \in \Z$ and $B^n, n \in \Z$ respectively. Then both 
${\cal A}$ and ${\cal B}$ are compact abelian groups. Let ${\cal C}$
denote the closure of $(A^{-1},B) ^n , n \in \Z$ in ${\cal A} \times
{\cal B}$.  Consider the  extension $\alpha$ of $a$ on $M \times {\cal
C} $  given by
\[ \alpha (x,y,z) = (a\, x, A ^{-1}\, y, B\, z).\] Since $(A^{-1},B)$
generates  a dense subgroup of ${\cal C}$, translation on
${\cal C}$ by $(A^{-1},B)$ is ergodic w.r.t. Haar measure. Since $a$
is weak  mixing, it follows that the extension $\alpha$ is still
ergodic w.r.t. the product of  Haar measure on ${\cal C}$ with the 
measure on $M$.   As above, let $M_{\varepsilon} \subset M$ be a set
of measure $1 -\varepsilon$ on which $\psi$ is continuous.  Since the
transformation induced by $\alpha$ on 
 $M_{\varepsilon} \times {\cal C}$ is ergodic, there is a point 
 $(x,y,z) \in M_{\varepsilon} \times {\cal C}$ whose $\alpha$-orbit
 intersected with $M_{\varepsilon} \times {\cal C}$ is dense
  in $M_{\varepsilon} \times {\cal C}$.
 Thus, given $x^* \in  M_{\varepsilon}$, there is  a sequence of
integers
 $n_k$ such that $a^{n_k} \, x \rightarrow x^*$, $a^{n_k} \, x \in
M_{\varepsilon}$,
 $A^{- n_k}\, y \rightarrow y$
 and $B^{n_k} \, z \rightarrow z$. It follows that  $A^{- n_k}
\rightarrow \mbox{Id}$  and 
 $B^{n_k} \rightarrow \mbox{Id}$. Since $\psi$ is continuous on 
$M_{\varepsilon}$ and $A$ commutes with the image of $\bar{\kappa}$,
equation $(*)$ implies that
\[ \psi (x^*) = \lim _{n_k \rightarrow \infty} \psi (a^{n_k} \, x)=
 \lim _{n_k \rightarrow \infty} B^{n_k} \, \psi (x) \, A^{-n_k} \,
\bar{\kappa}
 (a^{n_k},x)^{-1} = \psi (x) \, \lim _{n_k \rightarrow \infty}
 \bar{\kappa}(a^{n_k},x)^{-1} .\]

 Let $C$ be the compact subgroup of $H \subset GL(N,R)$ which contains
the image of $\bar{\kappa}$ and let 
$pr: {\cal G}
\rightarrow {\cal G}/C$ denote the projection . Then $pr \circ \psi $
is  constant a.e. on $M_{\varepsilon}$. 
 Since the $M_{\varepsilon}$ have measure arbitrarily close to 1, 
 $pr \circ \psi $ is constant a.e. on $M$.  Thus there is $\psi _0 \in
 {\cal G}$ such that for a.e. $x$, $\psi (x) \in \psi _0 \, C$.  By
the construction of $M'$ as $\frac{P \times H^0 \setminus H}{H}$,  the
$c$ frames $\bar{\tau} (x)=\tau (x) \, \psi _0$ belong to the fiber of
$P$ over $x \in M$. 
 Hence this fiber is just the $H$-orbit of any of
 the $c$ frames in $\bar{\tau} (x)$. In particular, we can choose the
 $H$-reduction $P$ smoothly. Hence $M'$ is a smooth finite manifold
cover
 of $M$ with a smooth $\Gamma$-action on $M'$. Pick a connected
component
 $M' _0$ of $M'$. Then some  subgroup $\Gamma _0$ of finite index 
 of $\Gamma$ preserves  $M' _0$. Passing to a finite power of $a$ if
necessary,
 we may assume that $a \in \Gamma _0$, and hence its lift to $M'$
preserves 
  $M'_0$. Let $\hat{a}$ denote this restriction. By choice of $a$,
$\hat{a}$ 
 is weak mixing  and also an affine automorphism of the homogeneous 
 space $M'_0$. The superrigidity framing on $M'$ restricts to a
superrigidity
 framing for $\Gamma _0$ on $M'_0$. Our previous argument for the
multivalued
 function $\psi: M \rightarrow {\cal G}$ can be applied to the single
valued function
 $M' _0 \rightarrow GL(N,\R)$ determined by a superrigidity framing.
Thus we get a constant translation $\tau '$ of the 
 standard frame on $M'_0$ which differs from a superrigidity framing
for the 
 $\Gamma _0$-action on $M'_0$ by translations by elements of $C$. Thus
$\tau '$
 is a superrigidity framing for the $\Gamma _0$-action on $M'_0$, a 
 conjugate of the representation $\pi$ and a cocycle taking values in
$C$.
\vspace{1em}
\QED

Now we are ready to prove the main theorem of this section. 

\begin{Theorem}   \label{thm-lattice-anosov}
 Let $G$ be  a linear semisimple  Lie 
group all of whose simple factors have real rank at 
least 2. Let $\Gamma$ be an irreducible lattice in $G$. Then a 
sufficiently small $C^1$-perturbation of an  algebraic  
Anosov action of $\Gamma$ on a nilmanifold $M$ is \ci -conjugate 
to the original action,  by a conjugacy $C^1$-close to the identity. 
\end{Theorem}

\def\N{{\cal N}}

\PROOF  Let $\rho$ be an algebraic Anosov action of $\Gamma$ on a nilmanifold
 $M$ of dimension $N$.  Thus its universal cover $\N$ is a nilpotent Lie group.

{\bf Step 1: Translation invariance of the superrigidity framing}
We will show that the hypotheses of Proposition~\ref{prop-locrig} hold.

Let us first show that the restriction of the action to $\Delta$ is 
still  Anosov. Fix a translation invariant  framing $\tau$ on $M$. 
Let $\sigma : \Gamma \rightarrow SL(N,\R)$ denote the
linear part of the $\Gamma$-action, and let $\tilde{\sigma}$ denote the
associated map into the adjoint group $PSL(N,\R)$.  Then $\rho (g)$ is an 
Anosov diffeomorphism precisely when $\tilde{\sigma }(g)$ does not have 
eigenvalues on the unit circle.  As $\rho$ is Anosov, this implies that 
$\tilde{\sigma} (\Gamma)$ is not a relatively compact subgroup of 
$PSL(N,\R)$. By Margulis' superrigidity theorem, $\tilde{\sigma} $ extends 
to a homomorphism $G \rightarrow PSL(N,\R)$ \cite[Theorem 6.16]{Margulis}. 

Now suppose that $g \in \Gamma$ is an Anosov element. Let $g_s$ be 
its semisimple Jordan component. Further decompose $g_s$ as a product
 $g_s = k\, p$  of commuting elements $k$ and $p$ where all eigenvalues of 
 $k$ have modulus 1, and all eigenvalues of the {\em polar part} $p$ are 
 positive. Since $G$ is an algebraic group, both $g_s$ and $p$ also belong to
 $G$. Since $\rho (g)$ is Anosov, $\tilde{\sigma} (g), \tilde{\sigma} (g_s)$ 
 and also $ \tilde{\sigma} (p)$ do not have 
eigenvalues on the unit circle. Let $D$ be a Cartan subgroup of $G$ 
containing $\Delta$. Let $A$ be the $\R$-split factor of $D$. Then $p$ is 
conjugate to some element $g' \in A$. 
Pick   $\delta$ in $\Delta$ such that the polar part of $\delta$ makes a 
sufficiently small angle with $g'$. Then $\sigma (\delta)$ also does not 
have eigenvalues on the unit circle. Hence $\rho (\delta)$ is Anosov.

Since some element  $a \in \Delta$ is Anosov, $\rho (a ) ^k$ for all $k$
is also Anosov and thus weak mixing on all connected finite covers of $M$.
It remains to prove   the local rigidity of
the action for some  Cartan subgroup $\Delta$ of 
$\Gamma$. We will
indicate two arguments for proving that. One argument works in general 
and uses  
Corollary~\ref{cor-nilautos}. 

Semisimplicity of the linear part follows from semisimplicity of $G$
 and 
the fact that $\Delta$ is a Cartan subgroup. The only remaining
condition is
 the absence of elements with roots of unity as
eigenvalues in the abelianization. This easily follows from
the following unpublished algebraic result of G. Prasad and A. 
 Rapinchuk  which strengthens a result of V. E. Voskresenskii \cite{Voskresenskii}.

\bigskip

\noindent {\bf Theorem} {\bf (Prasad-Rapinchuk)}   {\em Let $G$ and
$\Gamma$ be as before. Let  $\mu : G \rightarrow GL(N,\R)$ be a linear 
representation  of $G$. Then there exists a Cartan subgroup $\Delta$
of $\Gamma$ such that for all $\delta \neq 1 \in \Delta$,  the sum of
the multiplicities of roots of unity which are eigenvalues of 
$\mu (\delta)$ is the dimension of the 0-weight space of $\mu$. }
\bigskip

Let $\mu$ be the representation of $\Gamma$ on the abelianization 
$\frac{\frak h}{[\frak h, \frak h]}$ of the Lie algebra ${\frak h}$ of
$N$. Note
 that the kernel of
$\mu$ is finite by Margulis' finiteness theorem \cite[Theorem
4']{Margulis}. Hence $\pi$ is essentially faithful. By the last
 proposition and the presence of an Anosov element, we can pick a
Cartan  $\Delta 
\subset \Gamma$ such that $\mu (\lambda)$ does not have any roots of
unity 
 as eigenvalues for any $\lambda \in \Delta \neq 1$. Then 
 Corollary~\ref{cor-nilautos} applies
 to give a \ci -conjugacy between $\rho \mid _{\Delta} $ and the
 perturbation $\rho ^* \mid _{\Delta}$. 

If $\Gamma$ is a cocompact lattice in $G$, there is an alternate argument.
Suppose $\Gamma$ acts  on a nilmanifold $M$ by $\rho$, and that $\tilde{\rho}$
is a $C^1$-small perturbation of $\rho$. Induce $\rho$ and $\tilde{\rho}$ to
$G$ actions on $G \times _{\Gamma} M$. Restrict these actions to a maximal
split Cartan $A$. Let $K$ denote the compact part of the centralizer of $A$ in
$G$. Then both actions descend to $C^1$-close actions $\sigma$ and
$\tilde{\sigma}$ on  $K \setminus G \times _{\Gamma} M$. In fact, $\sigma$ is a
twisted Weyl chamber flow. By Corollary 5, there is a conjugacy $\phi$ close to
the identity and an automorphism $\alpha$ of $A$ close to the identity such
that
$\sigma = \phi \circ \tilde{\sigma} \circ \alpha$. Note that the fiber $M$ over
$K \, 1\, \Gamma$ is invariant under $\Delta := A \cap \Gamma$ (where we assume
that $\Gamma$ intersects $A$ in a lattice after a conjugation and moving the
base point  if necessary). Hence $\phi (M)$ is invariant under $\alpha ^{-1}
(\Delta)$. Note that $\phi (M)$ is close to $M$. Hence the projection $X$ of
$\phi (M)$ to $K \setminus G / \Gamma$ is an $\alpha ^{-1} (\Delta)$-invariant
compact set contained in an $\varepsilon$-neighborhood of $K\, 1 \, \Gamma$. By
the standard hyperbolic arguments, this set has to be contained in the neutral
manifold of the $A$ action on $K \setminus G / \Gamma$. If $\alpha$ is not the
identity, there is  $\delta \in \Delta$ such that $\alpha ^{-1} (\delta)  \,
\delta ^{-1}$ moves $K \, 1 \, \Gamma$ by more than $2 \, \varepsilon$. Then
$\alpha ^{-1} (\delta)$ has to move
$X$ outside the $\varepsilon$-neighborhood of of $K\, 1 \, \Gamma$ which is
impossible.

Thus $\alpha =id$ and $\phi$ is a conjugacy between $\tilde{\sigma} $ and 
a$\sigma$. Since $X$ is contained in the $A$-orbit of $K\, 1\, \Gamma$,
$\phi$ is a conjugacy of the suspension of the $\rho$ and 
$\tilde{\rho}$-actions to $\Rk$-actions. This implies that $\rho$ and 
$\tilde{\rho}$ are smoothly conjugate.

Thus all the conditions of Proposition~\ref{prop-locrig} have been verified, and
we get a translation invariant superrigidity framing $\tau ^*$ for a subgroup $\Gamma _0$ 
of finite index in $\Gamma$. We will also assume henceforth that $\rho$ and
$\tilde{\rho}$ agree on $\Delta$. Since $\tau ^*$ is a constant translate of
$\tau$ by some element $d \in GL(N,\R)$, we may further assume that 
$\tau = \tau ^*$, conjugating the relevant superrigidity representation and
cocycle by $d$, if necessary. 
\vspace{1em}

{\bf Step 2: Local \ci rigidity on a subgroup of finite index } 

By a result of S.  Hurder, there is a subgroup  $\Gamma ^*$ of finite index in 
$\Gamma$ which fixes a point $o$ in $M'$ under the unperturbed action $\rho$ 
\cite[Corollary 2]{Hurder2}.
G. A. Margulis showed that the first cohomology of $\Gamma^*$
with coefficients in any finite dimensional representation of $\Gamma^*$ 
vanishes \cite[Theorem 3']{Margulis}. Hence, by a theorem of D. 
Stowe, a sufficiently close perturbed action $\tilde{\rho}$ still has  a 
fixed point $o'$ for $\Gamma ^*$ \cite{Stowe}.
Since $\rho$ and $\tilde{\rho}$ coincide on $\Delta$, it follows that 
$\Delta \cap \Gamma^*$ fixes both $o$ and $o'$. Since $\Delta \cap \Gamma ^+$
contains Anosov elements and $o$ and $o'$ are close, this forces $o=o'$.
Hence we may assume that $\Gamma^*$ fixes $o$.

Consider the normal subgroup $\Gamma _A \subset \Gamma ^*$ which is 
generated by Anosov elements of $\Gamma ^*$ with respect to $\rho$. By 
Margulis' finiteness theorem, $\Gamma _A$ has finite index in $\Gamma$
\cite[Ch. IX, Theorem 5.4]{Margulis}. Since $\Gamma _A$ has Kazhdan's
property (T),
$\Gamma _A$ has a finite set $F$ of generators such that $\rho
(\gamma)$ is Anosov for all $\gamma \in F$ (this follows easily from
the standard proof of finite generation of Kazhdan groups
\cite[Theorem 7.1.5]{Zimmer}). 

Write $M=\N / \Lambda$ as a quotient of the universal cover $\N$ of
$M$.  Fix $\gamma$ in $F$, and lift $\rho (\gamma)$ and
$\tilde{\rho} (\gamma)$ to diffeomorphisms  $g$ and $\tilde{g}$ on the 
 universal cover $\N$ of $M$ such that both $g$ and $\tilde{g}$
 fix a point $p \in\N$ covering the common fixed point $o \in M$. 
Then  $g$ and 
$\tilde{g}$ are $C^1$-close on $\N$. Denote the lift of $\tau$ to $\N$
again by
$\tau$. Then $g$ transform $\tau$ by $\sigma (\gamma)$. By Step 1, 
$\tilde{g}$  transforms $\tau (x)$ by $\pi (\gamma) \circ \kappa
(\gamma, x)$.  Since 
 $\kappa$ takes values in a compact group $C$ commuting with the image
of $\pi$, 
 the stable and unstable distributions of $\tilde{g}$ (with respect to
a right
 invariant metric on $\N$) are determined by the
 eigenspaces of $\pi (\gamma)$, and   hence are invariant under right
translation by
 $\N$. Thus the stable  foliation of $\tilde{g}$ (resp. $g$) is an
 orbit foliation
 of some subgroup $\tilde{L}$ (resp. $L$) of $\N$. 
 The following lemma is presumably well known. For completeness we
include an 
 elegant  proof suggested to us by G. Prasad. 
 
\begin{Lemma}     \label{lem-Haus} Let $\N$ be a simply connected
nilpotent group, and $H$, $L$ two closed connected subgroups of
$\N$. Suppose that $L$ is contained in a tubular neighborhood of $H$ 
of fixed size
$\alpha$ with respect to a right invariant Riemannian metric on $N$.
Then 
$L$ is a subset of $H$.
\end{Lemma}

\PROOF Recall that $\N$ can be embedded into the upper triangular
matrices. Thus $\N , H$ and $L$ are unipotent algebraic groups. 
 Then the homogeneous space $H \setminus \N$ is an affine variety
\cite[Corollary 6.9]{Borel2}. Since $L$ is a unipotent group, the
$L$-orbit of $H \cdot 1$ in  $H \setminus \N$ is closed by Kolchin's
theorem \cite[Proposition 4.10]{Borel2}.  Hence this orbit is
homeomorphic with $(H \cap L) \setminus L$. Now suppose 
$L \subset H \, B$ where $B \subset \N$ is  compact. Then the
$L$-orbit of $H \cdot 1$ in  $H \setminus \N$ is bounded, and thus
compact. Since 
$H$ is unipotent, again $(H \cap L) \setminus L$ is an affine variety,
and thus cannot be compact. \QED

We will now show that the stable foliations for $g$ and $\tilde{g}$
coincide, or equivalently that $H = L$. Let $i$ be the injectivity
radius of $M$. We will assume that $\tilde{\rho}$ is so close to
$\rho$ that the unique topological conjugacy $\phi _{\gamma}$ between
$\rho (\gamma)$ and $\tilde{\rho} (\gamma)$ is $i/100$-close to the
identity. Note that $\phi _{\gamma}$ takes stable manifolds of $\rho
(\gamma)$ to stable manifolds of $\tilde{\rho} (\gamma)$. Since the
fixed points of $\rho (\gamma)$ are isolated, we may further assume
that $\phi _{\gamma} (o) =o$. Let $\bar{\phi} _{\gamma}$ be the lift
of 
$\phi _{\gamma}$ to $\N$ such that $\bar{\phi} _{\gamma} (p) =p$. Then
the distance between $\bar{\phi} _{\gamma}$ and the identity is again
at most 
$i/100$. Moreover,  $\bar{\phi} _{\gamma}$ takes the stable manifold
$L \, p$ to
$\tilde{L} \, p$. Thus the Hausdorff distance between $L$ and
$\tilde{L}$ is at most $i/100$.  Lemma~\ref{lem-Haus} now shows that
$L= \tilde{L}$.

Thus we now know that $g$ and $\tilde{g}$ have the same stable and
also  unstable manifolds ${\cal W} ^s (x) = L \, x$ and ${\cal W} ^u
(x) $. As before,
${\cal W} ^u$ is the orbit foliation of a subgroup, say $L^u$ of $\N$. 
 Set $f= g^{-1} \, \tilde{g}$. Then $f$ preserves 
 ${\cal W} ^s $ and ${\cal W} ^u $.  Since $\rho (\gamma)$ and
$\tilde{\rho} (\gamma)$ are $C^1$-close, their  induced maps on
$\Lambda = \pi _1 (M)$ coincide. Hence $f$ commutes with the action of
$\Lambda$ on $N$, and fixes $p$ as well as all $p \, \lambda $ for
$\lambda \in \Lambda$. 

Now note that 
${\cal W} ^u (p) \cap {\cal W} ^s (p\, \lambda)$ consists of at most
one point. Otherwise, $L \cap L^u \neq 1$, and hence contains a
one-parameter subgroup as the exponential map of $\N$ is a
diffeomorphism. That however contradicts the transversality of  ${\cal
W} ^s (p) =L\, p$ and 
${\cal W} ^u (p)= L^u \, p$.

Since $f$ fixes all $p \, \lambda$, it follows that $f$ also fixes any
intersection point ${\cal W} ^s (p\, \lambda) \cap {\cal W} ^u (p)$.
Since 
$\rho (\gamma)$ is Anosov, $W^s (o)$ is dense in $M$. Hence the set of 
intersection points ${\cal W} ^s (p\, \lambda) \cap {\cal W} ^u (p)$
is dense in
 ${\cal W} ^u (p)$. Thus $f$ fixes ${\cal W} ^u (p)$. As ${\cal W} ^u
(p)$ is
 dense in $M$, $f= id$, and hence $\rho (\gamma) = \tilde{\rho}
(\gamma)$.

 \vspace{1em}

{\bf Step 3: Local \ci -rigidity}

By Step 2,  we may assume that $\rho$  
and     $\tilde{\rho}$ coincide on $\Gamma _A$.
Now consider the normal subgroup $\Gamma _c \subset \Gamma$ which is 
generated by all Anosov elements of $\Gamma$ with respect to $\rho$. By 
Margulis' finiteness theorem, $\Gamma _c$ has finite index in $\Gamma$. 
Let $\gamma \in \Gamma _c$ be an Anosov element for $\rho$. Then some 
finite power $\gamma ^k$ belongs to $\Gamma  _A$. Since both $\rho (\gamma)$ 
and $\tilde {\rho} (\gamma)$ are close and commute with  
$\rho (\gamma) ^k = \tilde {\rho} (\gamma) ^k$, $\rho (\gamma)$ and
 $\tilde {\rho} (\gamma)$ coincide since the centralizer of an Anosov 
 element is discrete by expansiveness. 

To finish the proof of the proposition, we will  show that $\rho$  and     
$\tilde{\rho}$ automatically coincide on $\Gamma$. 
Let $\tau$ be the projection of a right-invariant framing as above. 
Fix $f \in \Gamma$, and  set $\tau ^f (x) = Df (\tau (f^{-1} (x))$.
Since $\Gamma _c \subset \Gamma$ is normal, one easily sees that  
$\tau ^f $ transforms under $\Gamma _c$ by the representation
$\sigma $ conjugated by the element $\sigma (f)$
\[ Dn \tau ^f (x) = \tau ^f (n\,x) \, \sigma (f) ^{-1} \, \sigma (n) \,
\sigma (f).\]
Thus the transformation law for $\tau ^f$ under $\Gamma _c$ is completely
given by a homomorphism. In particular, the cocycle $\kappa$ taking values 
in some compact group is trivial. Applying the argument of Step 2 to $\tau ^f$, 
one sees that $\tau ^f$ is a constant translate of $\tau$. This means that the
derivative $Df$ in terms of the framing $\tau$ is constant. Thus we get
a homomorphism $\sigma ^{\#} :\Gamma \rightarrow GL(N,\R)$ where $N$
is the
 dimension of $M$. 

We will now finish the proof of the theorem by showing that 
$\sigma $ and $\sigma ^{\#} $ coincide on $\Gamma$. The following argument for
 this
is due to G. A. Margulis, and relies on his structure theorem for 
homomorphisms of arithmetic groups in higher rank semisimple groups 
into algebraic groups \cite[Ch. VIII,Theorem 3.12]{Margulis}. 
Suppose that $\Gamma$ is an arithmetic group defined over a number
field $K$.  By Margulis' theorem, there is a homomorphism $\nu :
\Gamma \rightarrow  GL(N,\R)$ with finite image and a rational
homomorphism $\phi$ defined on the 
$\Q$-group obtained from $G$ by restriction of scalars from $K$ to $\Q$ such 
that $\sigma ^{\#} (\gamma) = \nu (\gamma) \phi (\gamma)$. Since $\phi$ and 
$\sigma$ are rational and coincide on the intersection of the kernel of $\nu$ 
with $\Gamma _A$, $\phi$ and $\sigma$ coincide. Since $\sigma ^{\#}$ and 
$\sigma$ coincide on $\Gamma _c$, $\nu$ is trivial on $\Gamma _c$. Hence 
$\nu$ is defined on the finite group $\Gamma / \Gamma _c$. On the other hand, 
 $\rho ^*$ is  $C^1$ -close to  $\rho$. Hence $\nu$ is arbitrarily close 
 to the identity. This is impossible as the trivial representation is 
 isolated within the space of representations of a fixed finite group. 
\QED

\section{Lattice actions on boundaries}

In this section we will demonstrate the local $C^{\infty}$-rigidity of
projective actions of cocompact lattices in semisimple groups of the
noncompact type.

\begin{Theorem} Let $G$ be a connected semisimple Lie group with 
finite center and without compact factors. Suppose that the real rank
of $G$ is at least 2.  Let $\Gamma \subset G$  be a cocompact 
irreducible lattice and $P$ a parabolic subgroup of $G$. Then the
action of $\Gamma$ on $G/P$ by left translations is locally
$C^{\infty}$-rigid. 
\end{Theorem}

The main idea in our proof is to reduce this rigidity problem to the
transversal  rigidity of the orbit foliation of a suitable normally
hyperbolic action  on $\Gamma \setminus G$. We would like to thank E.
Ghys for suggesting this approach. In fact, Ghys used a similar
duality to obtain the smooth classification of boundary actions of
Fuchsian groups \cite{Ghys}. 
 
In the special case of a projective action of a cocompact lattice  in
$SL(n+1,\R)$ on the $n$-sphere
$S^n$, M. Kanai established a local rigidity result for small
perturbations in the $C^4$-topology for all $n \geq 21$
\cite{Kanai}. His method is completely different and relies on the
vanishing  of a certain  cohomology. 

We begin with a brief review of foliated bundles and the suspension
construction (cf. \cite{HT} and also \cite{Yue} for a more detailed
description). Let $M$ be a compact manifold with  cover
$\tilde{M}$ (let us emphasize that $\tilde{M}$ need not be the
universal cover). Let $\Gamma$ be the group of covering
transformations.  Given a manifold $V$ and a homomorphism $\rho: \Gamma
\rightarrow \mbox{ Diff}(V)$, we can form the manifold  
\[ E = (\tilde{M} \times V )/ \Gamma \] where $\Gamma $ acts on
$\tilde{M} \times V$ by the diagonal action 
$\gamma (p, v)= (\gamma \, p,\rho(\gamma) v)$. Then  $E$ becomes  a
fiber bundle over $M$ with fiber
$V$ by projecting to the first factor. Note that the manifolds 
$\tilde{M} \times \{v\} $, $v \in V$,  project to a foliation ${\cal
F}$  of $E$ which is transverse to the fibers and has complementary
dimension.  This is the so-called  {\em suspension construction}. Any
bundle with a transversal foliation ${\cal F}$ of complementary
dimension is called a {\em foliated bundle}. One can show that any
such foliated bundle arises via the suspension construction where the
homomorphism $\rho :
\Gamma \rightarrow \mbox{ Diff }(V)$ is the {\em holonomy
homomorphism} of the foliation defined as follows. Let $p \in M$, and
identify $V$ with the fiber $V_p$. If $c$ is a loop at $p$, then for
$v \in V$, let $c (v)$ be the beginning point of a horizontal lift of
$c$ with endpoint $v$, i.e. a lift which is tangent to  ${\cal F}$.
This defines a diffeomorphism of $V$ which only depends on the
homotopy class of $c$.

Suppose now that $V$ is compact and that $\rho '$ is another action of
 $\Gamma$ on $V$. If $\rho'$ is sufficiently $C^1$-close to $\rho$ (on
some 
 finite generating
 set of $\pi_1 (M)$), then  the bundle $E'$ obtained from $\rho '$ via
the
 suspension construction is diffeomorphic to $E$ via a bundle map.
Moreover,
 the pullback of the natural transverse  foliation ${\cal F}'$ from
$E'$ to $E$
 is $C^1$-close to ${\cal F}$. We will henceforth identify $E'$ with
$E$, and
 ${\cal F}'$ with its pull-back. Finally, $\rho$ and $\rho '$ are
conjugate by
 a $C^{\infty}$ diffeomorphism  of $V$ precisely when there exists a 
 $C^{\infty}$ diffeomorphism of $E$ which carries ${\cal F}$ to ${\cal
F}'$.

Let $G$ be a connected semisimple Lie group with  finite center and
without compact factors. Let
$\Gamma \subset G$  be an irreducible lattice. We will now exhibit the
action of $\Gamma$ on a  boundary $G/P$ of $G$, $P$ a parabolic
subgroup,  as a holonomy action of a suitable ``stable'' foliation. 
We refer to \cite{Warner,Zimmer} for details on parabolic subgroups
and boundaries of
$G$.

Let $P=L\, C\, U^+$ be a Langlands  decomposition of $P$ (unique up to
conjugacy).  Here  $U^+$ is the unipotent  radical of $P$, and   
$L\,C$ is a Levi subgroup of $P$, i.e. a maximal reductive subgroup
(i.e. a product of a semisimple and an abelian group). Moreover, 
$C$ is the intersection of the Levi subgroup $L\,C$ with  a suitable
split Cartan $A$, and $L$ is reductive and commutes  with $C$.  If $P$
is a minimal parabolic then $U^+$ is a maximal unipotent subgroup of
$G$, $C$ is a split Cartan, and $L$ is just the maximal compact
subgroup of the centralizer of $A$ in $G$.

Let $c \in C$ be a regular element of $C$, i.e. an element such that
all roots of $G$ which are not zero on the Lie algebra ${\frak c}$  of
$C$ are also not zero on $\log c$. We may further assume that $c$
lies  on the boundary of a positive Weyl chamber of  $A$. The Lie
algebra of
$U^+$  is spanned by the rootspaces of $G$ whose roots are positive  
on ${\frak c}$. Hence these are precisely the roots which  are 
positive  on $\log c$. Moreover, the Lie algebra of $P$ is the sum of
the Lie  algebra of the centralizer of $A$ and the rootspaces which
are nonnegative on $\log c$. Let $P$ act on $\Gamma \setminus G$ by
right translations. Then  the orbit foliation  of $P$  is precisely 
the weak stable foliation of $c$. Let $G=K\, A \, N^+$ be an Iwasawa 
decomposition of
$G$, and set  $M_P = K \cap P$. Then $M_P$ is the centralizer of $c$
in $K$. Note that  $ \Gamma
\setminus G/M_P \rightarrow
\Gamma \setminus G/ K$ is a bundle with fiber $K/M_P $, and that  the
projection ${\cal W}^+$ of the orbit foliation of $P$ on $\Gamma
\setminus G$ to $\Gamma \setminus G  /M_P$
 is a transverse foliation of complementary dimension since $K\, P
=G$. As $C$ commutes with $M_P$, the right action of $C$ descends from
$\Gamma
\setminus G$ to $ \Gamma \setminus G/M_P $ with ${\cal W}^+$ as the
weak unstable foliation  of $c$.

Let us next describe the holonomy representation of ${\cal F}$  with
respect to
$\Gamma$ and the cover $G/M_P $ of $\Gamma \setminus G/  M_P$. The
fiber
$K/ M_P $ of the bundle 
$ \Gamma \setminus G/M_P \rightarrow \Gamma \setminus G/K $ is
diffeomorphic  with 
$G/P $, again since $G=K \, P$. Henceforth, we will identify the fiber
with  $G/P $.

\begin{Lemma} The holonomy representation  $\rho$ of the ${\cal
W}^+$-foliation  on the fiber $G/P$ is given by right multiplication 
by $\gamma$ on $G/ P $  for $\gamma \in \Gamma$.
\end{Lemma}

\proof Consider the cover $G/K $ of $\Gamma \setminus G/K $ and the
corresponding fiber bundle $G/ M_P  \rightarrow  G/K$. The
$K/M_P $-fiber of $1\, M_P $ gets identified with that of 
$\gamma \,M_P $ by left multiplication by $\gamma $. 
 Then $\rho (\gamma)(k\,M_P ) $  is that $l\,M_P$ such that $l \, P
\supset  \gamma\, k\, M_P  $.  Identifying $K/M_P $ with $G/P $, this
simply means that
$l\, P  = \gamma \, k \, P$. 

One can also see this more geometrically by identifying the orbit
foliation of
$P$ with the weak unstable foliation of $c$. Identify the fiber
$K/M_P $ first with a suitable subspace of the unit tangent sphere of 
$1 \, K$, and second with a subspace of the sphere at infinity of the
globally symmetric space $G/K$ by projecting the unit tangent vector
$v$ to  the asymptote class $v(-\infty )$ of the geodesic ray
$-v$ defines. One can
 check that the image is the homogeneous space $G/P $. Again, $v$ gets
identified with 
$d \gamma (v) \in T _{\gamma \, K} G/K$.  As $d \gamma (v) (- \infty)
= \gamma
\, x$, the holonomy image $\rho (\gamma) (v)$ is that $w \in T_{1\, K}
G/K$ with $w (\infty )=
\gamma \, x$. Under our identification of $T_{1\, K} G/K$ with $G/P$,
this just means $\rho (\gamma) (x) = \gamma \, x$. 
\QED

{\em Proof of Theorem 17}\,: Let us first describe the weak stable
foliation 
${\cal W}^- $ of $c$ and its holonomy. Let $U^-$ be the unipotent subgroup of $G$ whose  Lie algebra
is spanned by the rootspaces of $G$ whose roots are negative   on $c$ or equivalently ${\frak c}$.
Let $P^- =L\, C\, U^+$ be the opposite parabolic subgroup. Then $K/ M_P$ can be identified with
$G/P^-$. Geometrically, viewing $K/ M_P$ with a subset of $T_1 G/K$, this identification is given by
sending a suitable unit tangent vector $v$ to $v (- \infty)$.  This yields a diffeomorphism $\theta:
G/P \rightarrow G/P^-$. Now the holonomy representation $\rho ^-$ of ${\cal W} ^-$ is given by $\rho
^- (\gamma) (\theta (x) = \theta (\rho(\gamma) ( x))$.

Now  a $C^1$-small perturbation $\tilde{\rho} $ of $\rho$ gives a $C^1$-small  perturbation
$\tilde{\rho} ^- $ of $\rho ^-$. These perturbations of the  holonomies define
$C^1$-small perturbations $\tilde {{\cal W}}^+$ and $\tilde {{\cal W}}^- $ of the foliations ${\cal
W}^+$ and ${\cal W}^- $. Note that the orbit foliation 
${\cal O}$ of the Levi subgroup $H= L\, C$ of $P$ on $\Gamma \setminus G$ is invariant under the
action of $M_P$, and hence descends to a foliation ${\cal R}$ on
 $\Gamma \setminus G / M_P$. Moreover, ${\cal R}$ is the leafwise intersection of ${\cal W}^+$ and
${\cal W}^- $, and this intersection 
 is transversal. Hence the intersection of 
 $\tilde {{\cal W}}^+$ and $\tilde {{\cal W}}^- $ is a $C^1$-small perturbation 
 $\tilde { {\cal R}}$ of ${\cal R}$.  By Corollary 7, there is a
smooth  orbit 
 equivalence $\Phi$ from $\tilde{{\cal R}}$ to ${\cal R}$  
 that  
 $\Phi$ is 
 $C^1$-close to the identity. 
 Hence the pull-back foliations $\Phi ^{-1}(\tilde {{\cal W}}^+)$ and 
 $\Phi ^{-1}(\tilde {{\cal W}}^-)$ are $C^1$-close to ${\cal W}^+$ and 
 ${\cal W}^- $, and are also saturated under the ${\cal R}$ foliation.
 In particular, they are invariant under the action of $c$. However, it is clear
 that ${\cal W}^+$ and ${\cal W}^- $ are the unique $c$-invariant foliations
 $C^1$-close to ${\cal W}^+$ and ${\cal W}^- $. Indeed, any   
 vector close to the tangent space of  ${\cal W}^- $ but not tangent 
 to ${\cal W}^- $ has a non-trivial strong unstable component, and eventually gets
 expanded by $c$ under  forward  iteration. Thus eventually,
 the unstable  component would dominate, and the foliation could not be
 $C^1$-close to ${\cal W}^- $. 
 
 We conclude that $\Phi$ is a $C^{\infty}$-orbit equivalence between $\tilde {{\cal W}}^+$ 
 and ${\cal W}^+$. Hence the holonomy actions $\tilde{\rho}$ and $\rho$ are 
 $C^{\infty}$-conjugate. 
 \QED
  
\def \spa{\vspace{-.5em}}

Department of Mathematics, The Pennsylvania State University, University Park, PA 16802 \\
{\em E-mail address}: katok$\_$a$@$math.psu.edu

 Department of  Mathematics, University of Michigan, Ann Arbor, MI 48103\\ {\em E-mail
address}: spatzier$@$math.lsa.umich.edu


{\small
\begin{thebibliography}{aaaa}

\bibitem{A1} D.V.Anosov, {\em Roughness of geodesic flows on compact manifolds of negative
curvature} (in Russian),Dokl.Akad.Nauk SSSR, {\bf 145} (1962), 707--709. 
\spa
\bibitem{A2} D.V.Anosov,{\em Geodesic Flows on closed Riemannian manifolds of negative
curvature}, Ppoc. Steklov Math. Inst {\bf 90} (1967). 
\spa
\bibitem{Borel2} A. Borel, {\em Linear Algebraic Groups}, 2nd ed., Springer Verlag, New
York 1991.
\spa
\bibitem{Brezin-Moore} J. Brezin and C. C. Moore, {\em Flows on homogeneous spaces: a new
look}, Amer. J. of Math. {\bf 103} (1981), 571--613.
\spa
\bibitem{BS0}    K. Burns and R. J. Spatzier, {\em On  topological Tits buildings and their
classification},
 Publ. Math. IHES {\bf 65} (1987), 1--35.
\spa
\bibitem{Chen} K.-T. Chen, {\em Equivalence and decomposition of vector fields about an
elementary  critical point}, Amer. J. of Math. {\bf 85} (1963), 639--722.
\spa
\bibitem{Ghys} E. Ghys,  {\em Rigidit\'{e} differentiables des groupes fuchsiens}, Publ.
Math. IHES {\bf 78} (1993), 163--185.
\spa
\bibitem{GS1} E. Goetze and R. J. Spatzier, {\em Smooth classification of Cartan actions
of higher rank semisimple Lie groups and their lattices}, preprint.
\spa
\bibitem{HPS}           M. Hirsch, C. Pugh and M. Shub,
 {\em Invariant manifolds}, Lecture Notes in Mathematics  {\bf 583}, Springer Verlag,
Berlin  1977.
\spa

\bibitem{HT} M. Hirsch and W. P. Thurston, {\em Foliated bundles, invariant measures and
flat manifolds}, Ann. of Math. {\bf 101} (1975), 369--390.
\spa
\bibitem{Hurder1} S. Hurder, {\em Rigidity for Anosov actions of higher rank lattices},
Ann. Math. {\bf 135} (1992), 361--410.
\spa
\bibitem{Hurder2} S. Hurder, {\em Affine Anosov actions}, Mich. Math. J. {\bf 40} (1993),
561--575.
\spa
\bibitem{Kanai}  M. Kanai, {\em A new approach to the rigidity of discrete group actions},
preprint 1995.
\spa
\bibitem{K1} A. Katok, {\em Hyperbolic measures for actions of higher-rank abelian
groups}, preprint 1996. 
\spa
\bibitem{K2} A. Katok, {\em Normal forms and invariant geometric structures on transverse
contracting foliations}, preprint 1996. 
\spa
\bibitem{K3} A. Katok, {\em Local properties of hyperbolic sets} (in Russian), Supplement
to the Russian translation of ``Differentiable Dynamics'' by Z.Nitecki, Mir, Moscow, 
1975, 214--232.
\spa
\bibitem {KH} A.Katok and B.Hasselblatt, {\em Introduction to the modern theory of
dynamical systems},  Cambridge University Press, 1995
\spa
\bibitem{KL1} A. Katok and J. Lewis, {\em Local rigidity for certain groups of toral
automorphisms}, Israel J. of Math. {\bf 75} (1991), 203--241.
\spa
\bibitem{KL2} A. Katok and J. Lewis, {\em Global rigidity for lattice actions on tori and
new examples of volume preserving actions}, Israel J. of Math. {\bf 93}  (1996), 253--281. 
\spa
\bibitem{KLZ} A. Katok,  J. Lewis and R. J. Zimmer,  {\em  Cocycle superrigidity and
rigidity for lattice actions on tori}, Topology {\bf 35 } (1996),  27--38.
\spa
\bibitem{KS1}  A. Katok and R. J. Spatzier, {\em First cohomology of Anosov actions of
higher rank abelian  groups and applications to rigidity}, Inst. Hautes \'{E}tudes  Sci.
Publ. Math. No. {\bf 79} (1994), 131--156. 
\spa
\bibitem{KS2}  A. Katok and R. J. Spatzier, {\em Subelliptic estimates of polynomial
differential operators and  applications to rigidity of abelian action}, Math. Res. Lett.
{\bf 1} (1994), no. 2, 193--202.
\spa
\bibitem{KS3}  A. Katok and R. J. Spatzier, {\em Invariant Measures  for Higher Rank
Hyperbolic Abelian Actions}, Erg. Th. and Dynam. Syst. {\bf 16} (1996), no. 4,  751--778.
\spa
\bibitem{KS4}  A. Katok and R. J. Spatzier, {\em Differential rigidity of hyperbolic
abelian actions}, preprint 1992. 
\spa
\bibitem{KS5} A. Katok and R. J. Spatzier, {\em Non-stationary normal forms and rigidity
of group actions}, Elect. Res. Ann. of the AMS, {\bf 2} (1996), no. 3, 124--133.
\spa
\bibitem{Leuzinger} E. Leuzinger, {\em  New geometric examples of Anosov actions}, Ann.
Global Anal. Geom. {\bf 12} (1994), 173--181.
\spa
\bibitem{Lewis} J. Lewis, {\em The algebraic hull of the derivative cocycle}, preprint
(1993).
\spa
\bibitem{Mane} R.Ma\~{n}\'{e},  {\em A proof of the $C^1$ stability conjecture}, Inst.
Hautes \'{E}tudes Sci. Publ.
 Math. No. {\bf 66} (1987), 161-210. 
\spa
\bibitem{Margulis} G. A. Margulis, {\em Discrete subgroups of semisimple Lie groups},
Springer Verlag, Berlin 1991.
\spa
\bibitem{Mont-Zipp}  D. Montgomery and L. Zippin, {\em Topological
 transformation groups}, Interscience Publishers, Inc., New York, 1955.
\spa
\bibitem{Parry} W. Parry, {\em Ergodic properties of affine transformations and flows on
nilmanifolds}, Amer. J. Math. {\bf 91} (1969), 757--771.
\spa
\bibitem{Prasad-Ragh} G. Prasad and M. Raghunathan, {\em Cartan subgroups and lattices in
semisimple groups}, Ann. of Math. {\bf 96}  (1972), 296--317. 
\spa
\bibitem{PS} C. Pugh and M. Shub, {\em Ergodicity of Anosov actions}, Invent. Math. {\bf
15} (1972), 1--23.
\spa
\bibitem{Q1} N. Qian, {\em Anosov automorphisms for nilmanifolds and rigidity of group
actions},
 Erg. Th. and Dynam.  Syst. {\bf 15} (1995), 341--359.
\spa
\bibitem{Q2} N. Qian, {\em Topological deformation rigidity of actions of higher rank
lattices},
 Math. Res. Lett. {\bf 1} (1994), 485--499.
 \spa
\bibitem{Q3} N. Qian, {\em Tangential flatness and global rigidity of higher rank lattice
actions}, preprint.
\spa
\bibitem{Q4} N. Qian, {\em Smooth conjugacy for Anosov diffeomorphisms and rigidity of 
Anosov actions of higher rank lattices}, preprint.
\spa
\bibitem{Qian-Yue} N. Qian and C. Yue, {\em Local rigidity of Anosov higher rank lattice
actions}, preprint 1996.
\spa
\bibitem{Qian-Zimmer} N. Qian and R. J. Zimmer, {\em Entropy rigidity
 for semisimple group actions}, preprint.
\spa
\bibitem{Ragh} M. Raghunathan, {\em Discrete subgroups of Lie groups}, Springer, New York
1972. 
\spa
\bibitem{Robbin} J. Robbin, {\em A structural stability theorem}, Ann.Math {\bf 94}
(1971), 447--493. 
\spa
\bibitem{Sacksteder} R.Sacksteder, {\em Semigroups of expanding maps}, Trans. Amer. Math.
Soc. {\bf 221}  (1976), 281--288.
\spa
\bibitem{Seydoux} G. Seydoux, {\em Rigidity of ergodic volume-preserving actions of 
semisimple groups of higher rank on compact manifolds}, Trans. Amer. Math. Soc. {\bf 345}
(1994), no. 2, 753--776.
\spa
\bibitem{Smale} S.Smale {\em Differentiable dynamical systems}, Bull. Amer. Math. Soc.
{\bf 73} (1967), 747--817.
\spa
\bibitem{Sternberg}  S. Sternberg, {\em Local contractions and a theorem of Poincar\'{e}},
Amer. J. of Math. {\bf 79} (1957), 809--824.
\spa
\bibitem{Stowe} D. Stowe, {\em The stationary set of a group action}, Proc. AMS {\bf 79}
(1980), 139--146.
\spa
\bibitem{Voskresenskii} V. E. Voskresenskii, {\em Maximal tori without effect in
semisimple algebraic groups}, (Russian) Mat. Zametki {\bf 44} (1988), no. 3, 309--318.
\spa
\bibitem{Warner} G. Warner, {\em Harmonic analysis on semi-simple Lie groups I}, Springer
Verlag, Berlin 1972.
\spa
\bibitem{Yue} C. Yue, {\em Smooth rigidity of rank 1 lattice actions on the sphere at
infinity}, Math. Res. Lett. {\bf 2} (1995), 327--338.   
\spa
\bibitem{Zimmer}  R. J. Zimmer, {\em Ergodic Theory and Semisimple Groups}, Boston,
Birkh\"{a}user, 1984.
\spa
\bibitem{Zimmer1}  R. J. Zimmer, {\em On the algebraic hull of an automorphism group of a
principal bundle},
 Comment. Math. Helv. {\bf 65} (1990), 375--387.
\end{thebibliography} }
\end{document}